\documentclass{article}
\usepackage{amsmath,amsthm,amssymb}

\numberwithin{equation}{section}

{\theoremstyle{definition}\newtheorem{definitiona}{Definition}[section]

\newtheorem{terminology}[definitiona]{Terminology}
\newtheorem{remark}[definitiona]{Remark}
\newtheorem{example}[definitiona]{Example}}
\newtheorem{proposition}[definitiona]{Proposition}
\newtheorem{lemma}[definitiona]{Lemma}
\newtheorem{theorem}[definitiona]{Theorem}
\newtheorem{corollary}[definitiona]{Corollary}

\newenvironment{definition}{\begin{definitiona}}{\mbox{} \hfill
    $\blacktriangle$ \end{definitiona}}

\newcommand{\M}{\operatorname{M}}
\newcommand{\eindeformule}{\mbox{}\hfill\raisebox{0.6cm}[0mm][0mm]{$\blacktriangle$}\vspace{-0.6cm}}
\newcommand{\einde}{\hfill $\blacktriangle$}
\newcommand{\na}{\circ}
\newcommand{\bol}{\bullet}
\newcommand{\dia}{\, \sharp \,}
\newcommand{\car}{\diamond}
\newcommand{\ster}{\ast}
\newcommand{\actsl}{\triangleright}
\newcommand{\actsr}{\triangleleft}
\newcommand{\weak}{^{- \, \sigma\text{-weak}}}
\newcommand{\subq}{_{\hspace{-1.5mm}q}\hspace{1.5mm}} 
\newcommand{\subqq}{_{\hspace{-0.7mm}q}\hspace{0.7mm}} 
\newcommand{\cB}{\mathcal{B}}
\newcommand{\R}{\mathbb{R}}
\newcommand{\cO}{\mathcal{O}}
\newcommand{\cE}{\mathcal{E}}
\newcommand{\Npsi}{\mathcal{N}_\psi}
\newcommand{\de}{\delta}
\newcommand{\sde}{\Delta}
\newcommand{\te}{\theta}
\newcommand{\ot}{\otimes}
\newcommand{\al}{\alpha}
\newcommand{\be}{\beta}
\newcommand{\GGt}{G/G_2}
\newcommand{\GoG}{G_1 \backslash G}
\newcommand{\recht}{\rightarrow}
\newcommand{\Op}{^{\text{\rm op}}}
\newcommand{\Yh}{\hat{Y}}
\newcommand{\Wh}{\hat{W}}
\newcommand{\om}{\omega}
\newcommand{\B}{\mathcal{L}}
\newcommand{\K}{\mathcal{K}}
\newcommand{\si}{\sigma}
\newcommand{\io}{\iota}
\newcommand{\Mh}{\hat{M}}
\newcommand{\deh}{\hat{\delta}}
\newcommand{\cL}{\mathcal{L}}
\newcommand{\cR}{\mathcal{R}}
\newcommand{\cst}{C$^*$}
\newcommand{\Sh}{\hat{S}}
\newcommand{\sla}{\lambda}
\newcommand{\cW}{\mathcal{W}}
\newcommand{\pih}{\hat{\pi}}
\newcommand{\pitil}{\tilde{\pi}}
\newcommand{\cV}{\mathcal{V}}
\newcommand{\etah}{\hat{\eta}}
\newcommand{\C}{\mathbb{C}}
\newcommand{\ga}{\gamma}
\newcommand{\full}{_{\text{\rm f}}}
\newcommand{\red}{_{\text{\rm r}}}
\newcommand{\uni}{_{\text{\rm u}}}
\newcommand{\Jh}{\hat{J}}
\newcommand{\Ah}{\hat{A}}
\newcommand{\unit}{^\star}
\newcommand{\cAu}{\cA\unit}
\newcommand{\cA}{\mathcal{A}}
\newcommand{\cC}{\mathcal{C}}
\newcommand{\Si}{\Sigma}
\newcommand{\Q}{\mathbb{Q}}
\newcommand{\Z}{\mathbb{Z}}
\newcommand{\cP}{\mathcal{P}}
\newcommand{\cU}{\mathcal{U}}

\begin{document}
\begin{center}
{\Large\bf Non-semi-regular quantum groups coming from number theory}

\bigskip

{\sc by Saad Baaj, Georges Skandalis and Stefaan
  Vaes}
\end{center}

{\footnotesize Laboratoire de Math{\'e}matiques Pures;
Universit{\'e} Blaise Pascal; B{\^a}timent de Math{\'e}matiques; F--63177 Aubi{\`e}re
Cedex (France) \\ e-mail: Saad.Baaj@math.univ-bpclermont.fr \vspace{0.3ex}\\
Institut de Math{\'e}matiques de Jussieu; Alg{\`e}bres d'Op{\'e}rateurs et
Repr{\'e}sentations; 175, rue du Chevaleret; F--75013 Paris (France) \\
e-mails: skandal@math.jussieu.fr and vaes@math.jussieu.fr}

\bigskip

\begin{abstract}\noindent
In this paper, we study \cst-algebraic quantum groups obtained
through the bicrossed product construction.  Examples using groups of
adeles are given and they provide the first examples
of locally compact quantum groups which are not semi-regular: the
crossed product of the quantum group acting on itself by
translations does not contain any compact operator. We describe all
corepresentations of these quantum groups and the associated universal
\cst-algebras. On the way, we provide several remarks on
\cst-algebraic properties of quantum groups and their actions.
\end{abstract}

\section{Introduction}

What is the quantum analogue of a locally compact group? Several
authors addressed this question and came out with various
answers. G.I. Kac, Kac-Vainerman and Enock-Schwartz gave a quite
satisfactory set of axioms in \cite{E-S,Kac1,KVai} defining what Enock-Schwartz
called \emph{Kac algebras}. In the 1980's, new examples, which
were certainly to be considered as quantum
groups, but were not Kac algebras, were constructed, see e.g.\ \cite{Wor1}. Many
efforts were then made to enlarge the definition of Kac
algebras.

Already in the Kac algebra setting, it was known that all the
information on the quantum group could be encoded in one object: this
was called a \emph{multiplicative unitary} in \cite{BS}. It was then
quite natural to try and make constructions with a
multiplicative unitary as a starting point.

In order to be able to construct \cst-algebras out of a
multiplicative unitary, some additional assumptions were
made: this multiplicative unitary was assumed to be \emph{regular}
in \cite{BS}, \emph{semi-regular} in \cite{B}, \emph{manageable} in
\cite{Wor},~... In fact, a
multiplicative unitary can really be quite singular: in
Remark~\ref{strange.mult.un} below, we give
an example of a multiplicative unitary which should
certainly not be considered as the quantum analogue of
a locally compact group.

In \cite{KV1,KV2}, a definition of a locally compact quantum group
along the lines of Kac algebras, but with a weaker set of
axioms, was given and was shown to lead to a manageable
multiplicative unitary. We believe that the definitions of \cite{KV1,KV2} can be considered as the final ones, at least from the
measure theoretic (von Neumann) point of view.

On the other hand, regularity and
semi-regularity are natural conditions and present some
additional features -- that will be discussed below. In
particular, regularity is very much connected with the
Takesaki-Takai duality.

All this leaves us with many questions: what
are the relations between the properties of \emph{semi-regularity}
and \emph{manageability}? We actually really
tried hard to prove that these properties are equivalent. All
previously known examples were locally compact quantum groups
whose multiplicative unitary was semi-regular.

The main result of this paper is the construction of a
locally compact quantum group whose multiplicative unitary
is not semi-regular. This is done using a construction which
goes back to G.I. Kac \cite{Kac} and was used by several authors
\cite[...]{BS,BS2,Maj2,Tak,VV,VV2}. Let $G$ be a locally compact group and
let $G_1$ and
$G_2$ be closed subgroups such that the map
$\theta:(x_1,x_2)\mapsto x_1x_2$ from $G_1\times G_2$ into $G$
is a measure class isomorphism. Associated to this
situation, there is a locally compact quantum group. We will
show that its associated multiplicative unitary is
semi-regular if and only if $\theta$ is a homeomorphism from
$G_1\times G_2$ onto an open subset of $G$. Moreover, we
will be able to identify all the associated operator algebras.

In fact, in the examples considered up to now, $G$ was a
real Lie group and one could easily see that the associated
multiplicative unitary is always semi-regular (by
differentiability considerations).

The examples that we consider here are just adelic analogues of
these Lie groups. In particular, we will consider the case
were $G$ is an adelic $ax+b$ group and $G_1,G_2$ are natural
subgroups.

\medskip Note that in the non-regular case, the image
of $\theta$ is a countable union of compact sets with
empty interior and the complement of the image of $\theta$ has measure $0$. Although
this is a relatively common phenomenon in topology, it was quite a
surprise to us to encounter it  in the locally compact group
case. It was taken for granted in \cite{BS2} that such a phenomenon
could not occur. Because of this, the main result in \cite{BS2} is
incorrect as stated. One has either to add the extra
assumption that the associated multiplicative unitary is
semi-regular, or to change the conclusion: the product map
$(g_1,g_2)\mapsto g_1g_2$ is not a homeomorphism from
$G_1\times G_2$ onto a dense open subset of $G$ but a
measure class isomorphism. This is discussed in
Subsection~\ref{subsec.pentagonal}.

\medskip We conclude the introduction by explaining the structure of
the paper. In Section~\ref{sec.prelim}, we recall the definition and
main properties of multiplicative unitaries and locally compact
quantum groups. In
Section~\ref{sec.bicros}, we introduce the most general definition of
a matched pair of locally compact groups. We describe the associated
locally compact quantum groups obtained through the bicrossed product
construction. We compute the different associated \cst-algebras, as
well as the corepresentations and the covariant representations. We
give necessary and sufficient conditions for (semi-)regularity. In
Section~\ref{sec.semireg}, we give several examples, providing the
first examples of locally compact quantum groups that are not
semi-regular. Finally, in Section~\ref{sec.semireg}, we give different
characterizations of (semi-)regularity and we consider some
\cst-algebraic properties of coactions.

\section{Preliminaries} \label{sec.prelim}
In this paper, all locally compact groups are supposed to be second
countable. 

When $X$ is a
subset of a Banach space, we denote by $[X]$ the norm closed linear
span of $X$. We denote by $\K(H)$ and $\B(H)$ the compact, resp.\ the
bounded operators on a Hilbert space $H$. We use $\Si$ to denote the
flip map from $H \ot K$ to $K \ot H$, when $H$ and $K$ are Hilbert spaces.

By $\ot$, we denote several types of tensor products: minimal tensor
products of \cst-algebras, von Neumann algebraic tensor products or
Hilbert space tensor products. There should be no confusion.

\subsection{Multiplicative unitaries and locally compact quantum
  groups}

Multiplicative unitaries were studied in \cite{BS}. We have the
following definition.

\begin{definitiona}
A unitary $W \in \B(H \ot H)$ is called a \emph{multiplicative
  unitary} if $W$ satisfies the pentagonal equation
$$W_{12} \; W_{13} \; W_{23} = W_{23} \; W_{12} \; .$$
\eindeformule
\end{definitiona}

We associate to any multiplicative unitary two natural algebras $S$
and $\Sh$ defined by
\begin{equation} \label{eq.SandSh}
S = [ (\om \ot \io)(W) \mid \om \in \B(H)_* ] \; , \quad \Sh = 
[ (\io \ot \om)(W) \mid \om \in \B(H)_* ] \; .
\end{equation}
In
general, the norm closed algebras $S$ and $\Sh$ need not be
\cst-algebras, but they are when the multiplicative unitary is
regular or semi-regular, see Definition~\ref{def.regular} and
the papers \cite{B,BS}, or when it is manageable, see \cite{Wor}.

In these cases, we can define
comultiplications $\de$ and $\deh$ on the \cst-algebras $S$ and $\Sh$ by the formulas
$$\de: S \recht \M(S \ot S) : \de(x) = W(x \ot 1)W^* \; , \quad \deh :
\Sh \recht \M(\Sh \ot \Sh) : \deh(y) = W^*(1 \ot y)W \; .$$

Multiplicative unitaries appear most naturally as the right (or left)
regular representation of a locally compact (l.c.) quantum group. The
theory of locally compact quantum groups is developed in
\cite{KV1,KV2}.

\begin{definitiona}
A pair $(M,\de)$ is called a (von Neumann algebraic) l.c.\ quantum group when
\begin{itemize}
\item $M$ is a von Neumann algebra and $\de : M \recht M \ot M$ is
a normal and unital $*$-homomorphism satisfying the coassociativity relation : 
$(\de \ot \io)\de = (\io \ot \de)\de$.
\item There exist normal semi-finite faithful weights $\varphi$ and $\psi$ on $M$ such that
\begin{itemize}
\item $\varphi$ is left invariant in the sense that $\varphi \bigl( (\om \ot
\io)\de(x) \bigr) = \varphi(x) \om(1)$ for all $x \in M^+$ with
$\varphi(x) < \infty$ and all $\om \in M_*^+$,
\item $\psi$ is right invariant in the sense that $\psi \bigl( (\io \ot
\om)\de(x) \bigr) = \psi(x) \om(1)$ for all $x \in M^+$ with $\psi(x)
< \infty$ and all $\om \in M_*^+$. \einde
\end{itemize}
\end{itemize}
\end{definitiona}

There is an equivalent \cst-algebraic approach to l.c.\ quantum
groups and the link is provided by the right (or left) regular
representation. So,
suppose that $(M,\de)$ is a l.c. quantum group with right invariant
weight $\psi$. Represent $M$ on the GNS-space $H$ of $\psi$ and
consider the subspace $\Npsi \subset M$ of square integrable elements:
$$\Npsi = \{ x \in M \mid \psi(x^*x) < \infty \} \; .$$ Denote by
$\Gamma: \Npsi \recht H$ the GNS-map. Then, we define a unitary $V \in
\B(H \ot H)$ by the formula
$$V \bigl(\Gamma(x) \ot \Gamma(y) \bigr) = (\Gamma \ot \Gamma)\bigl(
\de(x) (1 \ot y) \bigr) \quad\text{for all}\quad x,y \in \Npsi \; .$$
The unitary $V$ is a multiplicative unitary and it is called the
\emph{right regular representation} of $(M,\de)$. The comultiplication
$\de$ is implemented by $V$ as above: $\de(x) = V(x \ot 1)V^*$.

Although $V$ need not be regular or semi-regular in general (see the
discussion below), it is always manageable (see \cite{KV1}) and we
have \cst-algebras $S$ and $\Sh$ defined by
Equation~\eqref{eq.SandSh} and we have the comultiplications $\de$ and $\deh$
on the \cst-algebraic level.

To compare notations between this paper and the papers \cite{KV1,KV2},
we observe that in \cite{KV1,KV2}, the left regular representation is
the main object. The $(S,\de)$ agrees with the $(A,\Delta)$ of \cite{KV1},
but $\Sh$ agrees with $\Jh \Ah \Jh$ and $\deh(y)$ agrees with $(\Jh
\ot \Jh) \hat{\Delta}(\Jh y \Jh) (\Jh \ot \Jh)$.

\subsection{Representations and corepresentations}

Next, we recall the notion of a representation and a corepresentation
of a multiplicative unitary, see Definition~A.1 in \cite{BS}.

\begin{definition}
Suppose that $W$ is a multiplicative unitary on the Hilbert space
$H$. We call $x \in \B(K \ot H)$ a representation of $W$ on the
Hilbert space $K$ if $x_{12} \; x_{13} \; W_{23} = W_{23} \; x_{12}$. We call
$y \in \B(H \ot K)$ a corepresentation on the Hilbert space $K$ if
$W_{12} \; y_{13} \; y_{23} = y_{23} \; W_{12}$.
\end{definition}

We remark that, if $W$ is the (left or right) regular representation
of a l.c.\ quantum group, then a representation $x$ satisfies $x \in \M(\K(K) \ot S)$ and a
corepresentation $y$ satisfies $y \in \M(\Sh \ot \K(K))$. The defining
relations become $(\io \ot \de)(x) = x_{12} \; x_{13}$ and $(\deh \ot
\io)(y) = y_{13} \; y_{23}$. It is clear that, conversely, if $x$ and $y$
satisfy these relations, they give a representation, resp.\ a
corepresentation of $W$.

In the same way as $C^*(G)$, we can define the universal \cst-algebras
$S\uni$ and $\Sh\uni$ for any l.c.\ quantum group (see \cite{JK}). There exists a
universal corepresentation $W\uni \in \M(\Sh \ot S\uni)$ such that $$S\uni =
[(\om \ot \io)(W\uni) \mid \om \in \B(H)_* ] \; . $$
There is a
bijective correspondence between representations $\pi : S\uni \recht
\B(K)$ of the \cst-algebra $S\uni$ and corepresentations $y$ of $W$ on
$K$ given by $y = (\io \ot \pi)(W\uni)$. All this is developed in \cite{JK}.

When we want to stress the distinction between the reduced and the
universal \cst-algebras, we denote $S$ by $S_{\text{\rm r}}$.

The following definition is taken from page 482 in \cite{BS}.
\begin{definitiona} \label{covariant}
Let $W$ be a multiplicative unitary. A pair $(x,y)$ of a
representation $x$ and a corepresentation $y$ of $W$ on the same
Hilbert space $K$ is called
covariant if $$y_{12} \; W_{13} \; x_{23} = x_{23} \; y_{12} \; .$$
\eindeformule
\end{definitiona}
We remark that Proposition~A.10 of \cite{BS} remains valid in the
setting of l.c.\ quantum groups. Hence, if $W$ is the (left or right)
regular representation of a l.c.\ quantum group and $(x,y)$ is a
covariant pair, then $x$ is stably isomorphic to $W$ as well as $y$,
but not jointly, and that is a very crucial point.

\subsection{Regularity and semi-regularity}

Let $V$ be a multiplicative unitary. Then, we have a naturally
associated algebra (see \cite{BS}), defined by
$$\cC(V) = \{ (\io \ot \om)(\Si V) \mid \om \in \B(H)_* \} \;.$$
We have $[\cC(V) \cC(V)] = [\cC(V)]$, but $[\cC(V)]$ is in general not
a \cst-algebra.

Next, we recall the notions of regularity \cite{BS} and
semi-regularity \cite{B} of a multiplicative unitary.

\begin{definition} \label{def.regular}
A multiplicative unitary $V$ is called regular if the closure of
$\cC(V)$ is $\K(H)$ and semi-regular if this closure contains $\K(H)$.
\end{definition}

When $V$ is the regular representation of a l.c.\ quantum group, we
have the following characterization of regularity and semi-regularity.

\begin{proposition} \label{char}
Let $V$ be the right regular representation of a l.c.\
quantum group on the Hilbert space $H$. Then, $[\cC(V)] \cong S
\rtimes\red \Sh$.

So,
$V$ is regular if and
only if $S \rtimes\red \Sh = \K(H)$ and $V$ is semi-regular if and
only if $S \rtimes\red \Sh$ contains $\K(H)$.
\end{proposition}

Remark that we consider $\de$ as a right coaction of $(S,\de)$ on the
\cst-algebra $S$. Then, the reduced crossed product $S \rtimes\red
\Sh$ is, by definition, given by $[\de(S) (1 \ot \Sh)]$. Using the
left regular representation, this last \cst-algebra is isomorphic to
$[S \Sh] \subset \B(H)$.

\begin{proof}
Suppose that $V$ is the right regular representation of $(M,\de)$.
From \cite{KV1,KV2}, we know that we have two modular conjugations $J$
and $\Jh$ at our disposal: $J$ is the modular conjugation of
the left invariant weight $\varphi$ and, up to a scalar, also the modular
conjugation of the right invariant weight $\psi$, while $\Jh$ is the
modular conjugation of the left invariant weight on the dual. We put
$U=J \Jh$. Then, $U$ is a unitary and $U^2$ is scalar. From
\cite{KV1,KV2}, we know enough formulas to apply Proposition~6.9 of
\cite{BS} and to conclude that $(\Si (1 \ot U) V)^3$ is
scalar. Hence, up to a scalar, we have
\begin{equation} \label{eq.irred}
\Si V = (U \ot 1) V^* (1 \ot U) \Si V^* \Si (U \ot 1) \; .
\end{equation}
Using the fact that $(J \ot \Jh)V(J \ot \Jh) = V^*$ and the fact that
$[(y \ot 1) \deh(x) \mid x,y \in \Sh] = \Sh \ot \Sh$, we get
\begin{align*}
[\cC(V)] &= [(\io \ot \om) \bigl( (1 \ot JyJ) \Si V (1 \ot \Jh x \Jh)
\bigr) \mid x,y \in \Sh, \om \in \B(H)_*] \\
&= [(\io \ot \om) \bigl( (\Si (J \ot \Jh) (y \ot 1) \deh(x) V^* (J \ot
\Jh) \bigr) \mid x,y \in \Sh, \om \in \B(H)_* ] = [\Jh \Sh \Jh \cC(V) ]
\; .
\end{align*}
Observing that $J \Sh J = \Sh$, we get that $[U \cC(V) U] = [\Sh U
\cC(V) U]$. We combine this last equality with
Equation~\eqref{eq.irred} and the fact that $V \in \M(\Sh \ot \K(H))$ to conclude that
$$[U \cC(V) U] = [(\io \ot \om) \bigl( (x \ot k) V^* (1 \ot U) \Si V^*
\Si \bigr) \mid x \in \Sh, \om \in \B(H)_*, k \in \K(H) ] = [\Sh S] \;
.$$
\end{proof}

\begin{corollary}
Let $V$ be the (left or right) regular representation of a l.c.\
quantum group on the Hilbert space $H$. Then, $[\cC(V)]$ is always a
\cst-algebra.
\end{corollary}

\begin{remark}
The isomorphism $[\cC(V)] \cong S \rtimes\red \Sh$ is typical for the
\emph{regular} representation. Suppose that, on the Hilbert space $K$,
$(x,y)$ is a covariant
representation of the right regular representation $V$ of a l.c.\
quantum group (in the sense of Definition~\ref{covariant}
below). Because both $x$ and $y$ are necessarily amplifications of the
regular representation (individually, but not jointly, see the remark
after Definition~\ref{covariant}), this means
that we have representations $\pi$ of $S$ and $\pih$ of $\Sh$ on $K$
such that $x = (\pih \ot \io)(V)$, $y = (\io \ot \pi)(V)$ and
which are covariant for the action $\de$ of $(S,\de)$ on $S$ in the
sense that $$(\pi \ot \io)\de(a) = x (\pi(a) \ot 1) x^* \quad\text{for
  all}\; a \in S \; .$$
So, we have a representation of the full crossed product $S
\rtimes\full \Sh$, whose image is $[\pi(S) \pih(\Sh)]$.

On the other
hand, we can write $\cV = (\pih \ot \pi)(V)$, which is a multiplicative
unitary on $K$. Then, the image of $S \rtimes\full \Sh$ is equal to
$[S_\cV \Sh_\cV]$. We claim that, $[\cC(\cV)]$ is Morita
equivalent to $[\cC(V)] \cong S \rtimes\red \Sh$. More precisely,
$[\cC(\cV)]$ is a \cst-algebra and $[\cC(\cV)] \ot \K(K \ot H)$ is
spatially isomorphic to
$[\cC(V)] \ot \K(K \ot K)$. Hence, up to Morita equivalence, $[\cC(\cV)]$
is always the same, while $[S_\cV \Sh_\cV]$ changes when $\cV$ runs
through the covariant images of $V$. Here, we can also observe that,
if $V$ is semi-regular or regular, then the same holds for any
covariant image $\cV$ of $V$.

From the proofs of Proposition 3.19 and
Corollaire 3.20 of \cite{B}, we know that the multiplicative unitaries
$\cV_{14} \cV_{24}$ and $V_{36}$ on $K \ot K \ot H$ are unitarily equivalent. It
follows that $[\cC(V)] \ot \K(K \ot K) \cong [\cC(\cV_{13} \cV_{23})] \ot
\K(H)$. An easy calculation shows that $B:=[\cC(\cV_{13} \cV_{23})] =
[\cC(\cV) \ot \K(K)] \Si \cV$. From the isomorphism above, we know that
$B$ is a \cst-algebra. Because $[(1 \ot \K(K))B] = B$ and $B$ is a
\cst-algebra, we have $B = [B (1 \ot \K(K))] = [\cC(\cV) \cC(\cV)] \ot
\K(K) = [\cC(\cV)] \ot \K(K)$. Our claim is proven.
\end{remark}

We make a more detailed study of regularity and semi-regularity of
l.c.\ quantum groups in Section~\ref{sec.semireg}.

\section{Bicrossed products of locally compact groups} \label{sec.bicros}

\begin{definition} \label{def.matched}
We call $G_1,G_2$ a matched pair of l.c.\ groups, if there is a given
l.c.\ group $G$ such that $G_1$ and $G_2$ are closed subgroups of $G$
satisfying $G_1 \cap G_2 = \{e\}$ and such that $G_1G_2$
has complement of measure $0$ in $G$.
\end{definition}

Suppose, throughout this section, that we have fixed a matched pair
$G_1,G_2$ of closed subgroups of $G$. We always use $g,h,k$ to denote
elements in $G_1$, $s,t,r$ to denote elements in $G_2$ and
accordingly, $dg, ds, \ldots$ to denote integration on $G_1, G_2$ with
respect to a fixed {\it right} Haar measure. We use $x,y$ to denote
elements in $G$ and $dx$ for the corresponding integration. The right
Haar measure on $G$ is fixed such that the following proposition
holds. We denote by $\sde_1,\sde_2$ and $\sde$ the modular functions
on $G_1,G_2$ and $G$.

\begin{proposition} \label{density}
The map $$\te : G_1 \times G_2 \recht G : \te(g,s) = gs$$ preserves
sets of measure zero and satisfies
$$
\int F(x) \; dx = \iint F(gs) \sde_1(g) \sde(g)^{-1} \; dg \; ds 
= \iint F(sg) \sde_2(s) \sde(s)^{-1} \; dg \; ds \; ,
$$
for all positive Borel functions $F$ on $G$.
\end{proposition}
\begin{proof}
The proof of Lemma~4.10 in \cite{VV} can essentially be copied: if $K
\in C_c(G_1 \times G_2)$, we define the bounded, compactly supported,
Borel function $\tilde{K}$ on $G$ by the formula $\tilde{K}(gs) =
K(g,s) \sde_1(g)^{-1} \sde(g)$ and $\tilde{K}(x)=
0$ for $x \in G \setminus G_1 G_2$. Then, we define $I(K) = \int
\tilde{K}(x) \; dx$. The same calculation as in \cite{VV} shows that
$I$ is a right invariant integral. Because we supposed that $G$ is
second countable, $G_1 \times G_2$ is $\sigma$-compact and so, the
integral $I$ must be non-zero. 
\end{proof}

From the point of measure theory, we can not distinguish $G_1 \times
G_2$ and $G$. We define almost everywhere on $G$ the Borel functions
$p_1:G \recht G_1$ and $p_2 : G \recht G_2$ such that $$x = p_1(x) \;
p_2(x) \; .$$ Then, we can define $$\tau : L^\infty(G_1) \ot
L^\infty(G_2) \recht L^\infty(G_2) \ot L^\infty(G_1) : \tau(F)(s,g) =
F(p_1(sg),p_2(sg)) \; .$$
The $^*$-automorphism $\tau$ will be a (von Neumann algebraic version
of an) inversion of $(L^\infty(G_1),\de_1)$ and
$(L^\infty(G_2),\de_2)$ in the sense of
\cite{BS}, Definition~8.1. Denoting by $\si$ the flip map, it
means that $\tau \si$ is a matching of $(L^\infty(G_2),\de_2)$ and 
$(L^\infty(G_1),\de_1^{\text{\rm op}})$ with trivial cocycles, in the
sense of \cite{VV}, Definition~2.1:
\begin{equation} \label{eq.matching}
(\tau \ot \io)(\io \ot \tau)(\de_1 \ot \io) = (\io \ot \de_1) \tau
\quad\text{and}\quad (\io \ot \tau)(\tau \ot \io)(\io \ot \de_2) =
(\de_2 \ot \io)\tau \; .
\end{equation}
Referring to Section~2.2 and in
particular, Theorem~2.13 of \cite{VV}, we can construct a l.c.\
quantum group using $\tau$. Because we have chosen to follow
systematically the conventions of \cite{BS} and \cite{BS2}, we will
state explicitly the needed formulas.

Define $\al : L^\infty(G_1) \recht L^\infty(G_2) \ot L^\infty(G_1)$
and $\be : L^\infty(G_2)
\recht L^\infty(G_2) \ot L^\infty(G_1)$ by $$\al(F) = \tau(F \ot 1) \quad\text{and}\quad  \be(F) = \tau(1 \ot F) \; .$$
Then, $\al$ defines a left action of $G_2$ on $L^\infty(G_1)$ while
$\be$ defines a right action of $G_1$ on $L^\infty(G_2)$. In fact, if
we consider $\tau \si$ as a matching of $(L^\infty(G_2),\de_2)$ and 
$(L^\infty(G_1),\de_1^{\text{\rm op}})$ in the sense of \cite{VV},
these two actions $\al$ and $\be$ precisely agree with the actions
$\al$ and $\be$ appearing in Definition~2.1 of \cite{VV}. We also have
the obvious formulas $$\al(F)(s,g) = F(p_1(sg)) \quad\text{and}\quad
\be(F)(s,g) = F(p_2(sg)) \; .$$
If we equip the quotient spaces $\GGt$
and $\GoG$ with their canonical invariant measure class, the embedding
$G_1 \mapsto \GGt$ identifies $G_1$ with a Borel subset of $\GGt$ with
complement of measure zero. Because of Proposition~\ref{density}, this
embedding, as well as its inverse, respects Borel sets of measure
zero. Hence, it induces an isomorphism between $L^\infty(G_1)$ and $L^\infty(\GGt)$.
Then, $p_1$ can be considered as the
projection of $G$ to $\GGt$. Any $s \in G_2$ acts as a homeomorphism
on $\GGt$, preserving Borel sets of measure zero and hence, the map
$p_1(s \; \cdot)$ is defined almost everywhere on $G_1$, preserves Borel
sets of measure zero and gives
precisely the automorphism $\al_s$ of $L^\infty(G_1)$. Similar
considerations can be made for $\be$.

We also conclude that, for a fixed $s \in G_2$, $sg \in G_1 G_2$ for
almost all $g \in G_1$. An analogous statement holds for a fixed $g
\in G_1$. It follows that $p_1(sg)$ and $p_2(sg)$ are defined almost
everywhere if either $s$ or $g$ is fixed. This will allow us to use
freely $p_1(sg)$ and $p_2(sg)$ in integrals, see e.g.\ Equation~\eqref{eq.voorbeeld}.

We denote by
$X$ and $Y$ the usual multiplicative unitaries on $H_1:=L^2(G_1)$ and
$H_2:=L^2(G_2)$ respectively, defined by $(X \xi)(g,h) = \xi(gh,h)$
and $(Y \eta)(s,t) = \eta(st,t)$. Then $\Yh$ is the multiplicative
unitary on $L^2(G_2)$ given by $(\Yh \eta)(s,t)  = \sde_2(s)^{1/2} \eta(s,s^{-1}t)$.

Following \cite{BS} and \cite{VV}, we can define the main actors
of the paper.

\begin{definition} 
Defining
$$W = (\be \ot \io)(\Yh)_{123} (\io \ot \al)(X)_{234} \; ,$$
we get a multiplicative unitary $W$, which is the right regular
representation of the l.c.\ quantum group $(M,\de)$ given by
$$M =  (\al(L^\infty(G_1)) \cup \cL(G_2) \ot
  1)^{\prime\prime} = G_2 \ltimes_\al L^\infty(G_1)
  \quad\text{with}\quad \de(z) = W (z \ot 1) W^* \; ,$$
where $\cL(G_2)$ denotes the left regular representation of the group
von Neumann algebra of $G_2$, generated by the unitaries $\lambda_s$
defined by $(\lambda_s \xi)(t) = \Delta_2(s)^{1/2} \xi(s^{-1} t)$.
\end{definition}

The right Haar measure of $G_1$ is the left Haar measure of $G_1\Op$
and secondly, using the modular function $\sde_2$, the $L^2$-spaces of $G_2$
with right or left Haar measure are naturally isomorphic. Under this
isomorphism $W$ agrees exactly with the $\Wh$ of \cite{VV},
Definition~2.2 (one just has to interchange the indices $1$ and $2$
referring to $G_1$ and $G_2$ everywhere). Hence, the $(M,\de)$ defined
above agrees with $(M,\Delta\Op)$ in \cite{VV}, where $\Delta\Op =
\sigma \Delta$ with $\sigma$ the flip map.
So, we get indeed that $W$ is the
right regular representation of the l.c.\ quantum group $(M,\de)$.

Dually, we have
$$\Mh =  (\be(L^\infty(G_2)) \cup 1 \ot
  \cR(G_1))^{\prime\prime} = L^\infty(G_2) \rtimes_\be G_1 \quad\text{with}\quad
\deh(z) = W^* (1 \ot z) W \; ,$$
where $\cR(G_1)$ denotes the right regular representation of the group
von Neumann algebra of $G_1$, generated by the unitaries $(\rho_g
\xi)(h) = \xi(hg)$.

Then, $(\Mh,\deh)$ above and $(\Mh,\hat{\Delta})$ in \cite{VV} really agree.

From Equation~\eqref{eq.matching}, we get the following formulas for
the comultiplication on the generators of $M$ and $\Mh$:
\begin{align}
& \de \al = (\al \ot \al) \de_1 \quad\text{and}\quad (\io \ot
\de)(\Yh_{12}) = \Yh_{12} \; ((\io \ot \al)\be \ot \io)(\Yh)_{1234} \; ,
\label{eq.generatorM} \\
& \deh \be = (\be \ot \be) \de_2 \quad\text{and}\quad (\deh \ot
\io)(X_{23}) = (\io \ot (\be \ot \io)\al)(X)_{2345} \; X_{45} \; .
\label{eq.generatorMhat}
\end{align}

\begin{remark} \label{interchange}
If we interchange $G_1$ and $G_2$ and keep $G$, we also have a matched
pair of l.c.\ groups. Performing the bicrossed product construction
again, we get a multiplicative unitary on $L^2(G_1 \times G_2)$, which
through the unitary $(v \xi)(s,g) = \Delta_1(g)^{1/2}
\Delta_2(s)^{1/2} \xi(g^{-1},s^{-1})$ corresponds to $\Si W^*
\Si$. So, $S$ and $\Sh$ are interchanged and the comultiplications are flipped.
\end{remark}

So, the von Neumann algebraic picture of the l.c.\ quantum group
$(M,\de)$ is completely clear. In this paper, we study the associated
\cst-algebras. To determine these, the following lemma will be crucial.

\begin{lemma} \label{quotient}
Using the embedding $C_0(\GGt) \recht \M(C_0(G_1)) \subset
L^\infty(G_1)$, we have
$$[(\om \ot \io)\al(F) \mid \om \in \B(H_2)_*, F \in C_0(G_1) ] =
C_0(\GGt) \; .$$
\end{lemma}
\begin{proof}
We have to prove that for all $F_i \in C_c(G_i)$, the function
\begin{equation} \label{eq.voorbeeld}
H(g) := \int F_1(p_1(sg)) F_2(s) \; ds
\end{equation}
belongs to $C_c(\GGt)$ and
that these functions $H$ span a dense subspace of $C_0(\GGt)$. It
suffices to look at $F_1 = \tilde{K}_1 * P_1$ with $K_1,P_1 \in
C_c(G_1)$, i.e.
$$F_1(g) = \int K_1(h) P_1(hg) \; dh \; .$$ Take some $P_2 \in
C_c(G_2)$ such that $\int P_2(s^{-1}) \; ds = 1$. On $G$, we define bounded
Borel functions $K$ and $P$ with compact support by $$K(hs) = K_1(h)
F_2(s) \sde_1(h)^{-1} \sde(h) \quad\text{and}\quad P(hs) = P_1(h)
P_2(s)$$ and such that $K$ and $P$ equal $0$ outside $G_1 G_2$. In
particular, $K,P \in L^2(G)$ and defining $Q = \tilde{K} * P$, i.e.\ $$Q(x) =
\int K(y) P(yx) \; dy \; ,$$ we have $Q \in C_c(G)$.

We claim that $$H(g) = \int Q(gt^{-1}) \; dt \; .$$ Take $g \in
G_1$. Then, using Proposition~\ref{density},
\begin{align*}
\int Q(gt^{-1}) \; dt &= \iint K(y) P(y g t^{-1}) \; dy \; dt \\ &=
\iiint K(hs) P(hsgt^{-1}) \sde_1(h) \sde(h)^{-1} \; dh \; ds \; dt \\
&= \iiint K_1(h) F_2(s) P_1(h p_1(sg)) P_2(p_2(sg) t^{-1}) \; dt \; dh
\; ds \\ &= \iint K_1(h) F_2(s) P_1(h p_1(sg)) \; dh \; ds = H(g) \; .
\end{align*}
This proves our claim. We conclude that $H \in C_c(\GGt)$. Moreover,
the functions $H$ span a dense subset of $C_0(\GGt)$, because in the
proof above, the functions $K$ and $P$ both span a dense subset of
$L^2(G)$, so that the functions $Q$ span a dense subset of $C_0(G)$.
\end{proof}

From \cite{VV}, we know that the multiplicative unitary $W$ is the
(right) regular representation of the l.c.\ quantum group $(M,\de)$. Hence, we have two
associated \cst-algebras:
\begin{equation} \label{cstdefin}
S = [ (\om \ot \io)(W) \mid \om \in \B(H_2 \ot H_1)_* ] \quad\text{and}\quad \Sh = [ (\io \ot \om)(W) \mid \om \in
\B(H_2 \ot H_1)_* ] \; .
\end{equation}
The comultiplications $\de$ and $\deh$ restrict nicely to morphisms
$\de : S \recht \M(S \ot S)$ and $\deh : \Sh \recht \M(\Sh \ot \Sh)$.

As a first application of Lemma~\ref{quotient}, we can identify
$S$ and $\Sh$. We use $\lambda$ to denote the left regular
representation of $C^*_r(G_2)$ on $L^2(G_2)$.

\begin{proposition} \label{cst}
The following equalities hold:
\begin{align*}
& S = G_2 \ltimes\red C_0(\GGt) = \bigl[ (\sla(C^*_r(G_2)) \ot 1)
\al(C_0(\GGt)) \bigr] \quad\text{and}\quad \al : C_0(G_1) \recht \M(S)
\quad\text{is non-degenerate} \; ,
\\
& \Sh = C_0(\GoG) \rtimes\red G_1 = \bigl[ (1 \ot C^*_r(G_1))
\be(C_0(\GoG)) \bigr] \quad\text{and}\quad \be : C_0(G_2) \recht \M(\Sh)
\quad\text{is non-degenerate} \; .
\end{align*}
\end{proposition}
\begin{proof}
Observe that $\be(F) = B(F \ot
1)B^*$ for all $F \in L^\infty(G_2)$ where $B \in \M(\K(H_2) \ot
C_0(G_1))$ is a representation of $X$, the canonical implementation of
$\be$. Because $B$ is a representation of $X$, we get
\begin{align*}
S & = [(\om \ot \mu \ot \io \ot \io)(\Yh_{13} (\io \ot \io \ot
\al)(B^*_{12} X_{23})) \mid \om \in \B(H_2)_*, \mu \in
\B(H_1)_*] \\ &= [(\om \ot \mu \ot \io \ot \io)(\Yh_{13} (\io \ot \io \ot
\al)(B_{13} X_{23})) \mid \om \in \B(H_2)_*, \mu \in
\B(H_1)_*] \\ &= [(\om \ot \io \ot \io)\bigl(\Yh_{12} (\io \ot \al)(B (1
\ot C_0(G_1))) \bigr) \mid \om \in \B(H_2)_*] = 
[(\sla(C^*_r(G_2)) \ot 1) \al(C_0(G_1))] \; .
\end{align*}
Because $S$ is a \cst-algebra and because of Lemma~\ref{quotient}, we
get
$$S = [(\sla(C^*_r(G_2)) \ot 1) \al(C_0(G_1)) (\sla(C^*_r(G_2)) \ot
1)] = [(\sla(C^*_r(G_2)) \ot 1) \al(C_0(G/G_2)) (\sla(C^*_r(G_2)) \ot
1)] \; .$$
By definition, $G_2 \ltimes\red C_0(\GGt) = [ (\sla(C^*_r(G_2)) \ot 1)
\al(C_0(\GGt))]$ and this is a \cst-algebra. Hence, $S = 
G_2 \ltimes\red C_0(\GGt)$. So, we have proven
the first line of the statement above. The
second line is analogous (or follows by interchanging $G_1$ and $G_2$
and applying Remark~\ref{interchange}).
\end{proof}

Next, we prove the corresponding result for the universal
\cst-algebras.

\begin{proposition} \label{universalcst}
Suppose that $\cW$ is a corepresentation of $W$ on a Hilbert space $K$.
\begin{enumerate}
\item There exist unique
corepresentations $y$ of $\Yh$ and  $x$ of $X$ on $K$ such that
$$\cW = (\be \ot \io)(y) \; x_{23} \; . $$ 
\item
If we denote by
$\pi_1$ and $\pih_2$ the representations of $C_0(G_1)$
and $C^*(G_2)$ corresponding to $x$ and $y$ and if
we define the representation $\pitil_1$ of $C_0(\GGt)$ using the
extension of $\pi_1$ to $\M(C_0(G_1))$, then the pair
$(\pitil_1,\pih_2)$ is a covariant representation for the left action
of $G_2$ on $C_0(\GGt)$.
\item
The representation $\pi$ of
$S\uni$ associated to $\cW$ and the representation $\pitil_1 \times \pih_2$ of $G_2 \ltimes\full C_0(\GGt)$ have
the same image. 
\end{enumerate}
Conversely, every covariant representation
for the left action of $G_2$ on $\GGt$ is obtained in this way from a
corepresentation $\cW$ of $W$.

In particular, $$S\uni \cong G_2 \ltimes\full C_0(\GGt)$$
in a natural way.
An analogous statement holds for $\Sh\uni$ and in particular,
$\Sh\uni \cong C_0(\GoG) \rtimes\full G_1 \; .$
\end{proposition}
\begin{proof} {\it a)} This follows from \cite{BV}. Observe that $x
  \in \M(C^*_r(G_1) \ot \K(K))$ and $y \in \M(C_0(G_2) \ot \K(K))$.

{\it b)} Because $\cW$ is a
corepresentation, we get the following equality in $\B(H_2 \ot H_1 \ot
H_2 \ot H_1 \ot K)$.
\begin{equation}\label{formula.2}
(\deh \ot \io)(x_{23}) = (\be \ot \io)(y^*)_{345} \; x_{25}\; 
(\be \ot \io)(y)_{345} \; x_{45} \; .
\end{equation}
Define the non-degenerate
representation $\pitil_1$ of $C_0(\GGt)$ as in the statement of the
proposition.

We know that, in $\B(H_2 \ot H_1 \ot
H_2 \ot H_1 \ot H_1)$,
\begin{equation} \label{formula.1}
(\deh \ot
\io)(X_{23}) = (\io \ot (\be \ot \io)\al)(X)_{2345} \; X_{45} \; .
\end{equation}
From Proposition~\ref{cst}, it follows that $X_{23} \in \M(\Sh \ot
C_0(G_1))$ and so,
it follows that $(\io \ot (\be \ot \io)\al)(X)_{2345} \in \M(\Sh \ot
\Sh \ot C_0(G_1))$, from which we conclude that
\begin{equation}  \label{formula.non-deg}
(\be \ot \io)\al : C_0(G_1) \recht \M(\Sh \ot C_0(G_1))
\end{equation}
is a
well-defined and non-degenerate $^*$-homomorphism. From
Equation~\eqref{formula.1} and the fact that $(\io \ot \pi_1)(X) = x$,
we conclude that, in $\B(H_2 \ot H_1 \ot
H_2 \ot H_1 \ot K)$,
$$(\deh \ot
\io)(x_{23}) = (\io \ot (\io \ot \io \ot \pi_1)(\be \ot
\io)\al)(X)_{2345} \; x_{45} \; .$$
Combining with Equation~\eqref{formula.2}, we find that
\begin{equation} \label{formula.8}
(\io \ot (\io \ot \io \ot \pi_1)(\be \ot
\io)\al)(X) = (\be \ot \io)(y^*)_{234} \; x_{14}\; 
(\be \ot \io)(y)_{234} \; ,
\end{equation}
which yields in $\B(H_2 \ot H_1 \ot K)$ the equality
\begin{equation} \label{formula.3}
(\io \ot \io \ot \pi_1)(\be \ot
\io)\al(F) = (\be \ot \io)(y^*) (1 \ot 1 \ot \pi_1(F)) (\be \ot
\io)(y)
\end{equation}
for all $F \in C_0(G_1)$. Because of
non-degenerateness, the same holds for all $F \in \M(C_0(G_1))$ and
hence, for $F \in C_0(\GGt)$. As $G_2$ acts continuously on
$C_0(\GGt)$, we have $\al : C_0(\GGt) \recht \M(C_0(G_2) \ot
C_0(\GGt))$ and so, we arrive at the formula
$$(\io \ot \pitil_1)\al(F) = y^* (1 \ot \pitil_1(F)) y \;
,$$ for all $F \in C_0(\GGt)$. This precisely gives the required
covariance of $(\pitil_1,\pih_2)$.

{\it c)}
Denote by $\pi$ and $\pitil_1 \times \pih_2$ the corresponding
representations of $S\uni$ and $G_2 \ltimes\full C_0(\GGt)$. We have to
prove that they have the same image and this will be analogous to the
proof of Proposition~\ref{cst}. We use again the canonical
implementation $B \in \M(\K(H_2) \ot
C_0(G_1))$ of $\be$.
Observe that $$\pi(S\uni) = [(\om \ot \mu \ot \io)(\cW)
\mid \om \in \B(H_2)_*, \mu \in \B(H_1)_* ] = [(\om \ot
\mu \ot \io)(y_{13} \; B^*_{12} \; x_{23})\mid \om \in \B(H_2)_*, \mu \in \B(H_1)_*] \; .$$
Because $B_{12} \; B_{13} \; X_{23} = X_{23} \; B_{12}$, we get $B^*_{12} \; x_{23}
= (\io \ot \pi_1)(B)_{13} \; x_{23} \; B^*_{12}$ and hence,
\begin{align*}
\pi(S\uni) &= [(\om \ot \mu \ot \io)(y_{13} \; (\io \ot \pi_1)(B)_{13} \; x_{23})
\mid \om \in \B(H_2)_*, \mu \in \B(H_1)_* ] \\ &= [(\om \ot \io)(y \; 
(\io \ot \pi_1)(B) \; \pi_1(C_0(G_1)) \mid \om \in \B(H_2)_*] =
[\pih_2(C^*(G_2)) \; \pi_1(C_0(G_1)) ] \; .
\end{align*}
From Equation~\eqref{formula.3}, it follows that
$$\pitil_1((\om \be \ot \io)\al(F)) = (\om \ot \io) \bigl(
(\be \ot \io)(y^*) \; (1 \ot 1 \ot \pi_1(F)) \; (\be \ot
\io)(y) \bigr)$$ for all $F \in C_0(G_1)$ and $\om \in \B(H_2 \ot
H_1)_*$. Combining this with Lemma~\ref{quotient}, it follows that
$$[\pih_2(C^*(G_2)) \; \pitil_1(C_0(\GGt)) \; \pih_2(C^*(G_2))] =
[\pih_2(C^*(G_2)) \; \pi_1(C_0(G_1)) \; \pih_2(C^*(G_2))] \; .$$
The image of $\pitil_1 \times
\pih_2$ is given by $[\pih_2(C^*(G_2)) \; \pi_1(C_0(\GGt))]$ and is a
\cst-algebra. Also $\pi(S\uni)$ is a \cst-algebra, so that the
calculation above shows that, indeed $\pi$ and $\pitil_1 \times
\pih_2$ have the same image. This concludes the proof of the first
part of the proposition.

Suppose now, conversely, that we have a covariant representation
$(\pitil_1,\pih_2)$ on $K$ for the left action of $G_2$ on $\GGt$. Denote by
$y$ the corepresentation of $\Yh$ corresponding to $\pih_2$. The
representation $\pitil_1$ defines a measure class on $\GGt$ which, by
covariance, is invariant under the action of $G_2$. We claim that this
measure class is supported by the image of $G_1$ in $\GGt$. By
invariance, the transformation $(s,\bar{x}) \mapsto (s,\overline{sx})$ of
$G_2 \times \GGt$ preserves sets of measure zero and so, we have to prove that
the set of pairs $(s,\bar{x})$ such that $\overline{sx} \not\in G_1$ has
measure zero. Using the Fubini theorem, it is enough to prove that, for any $x \in
G$, the set of all $s \in G_2$ such that $\overline{sx} \not\in G_1$ has
measure zero. But this last set is equal to the set of all $s \in G_2$
such that $gs \not\in G_1 G_2 x^{-1}$ for all $g \in G_1$. And here,
we get a set of measure zero and our claim has been proved.

Because the measure class on $\GGt$ is supported by $G_1$, we can use
the Borel calculus to find a non-degenerate representation $\pi_1$ of
$C_0(G_1)$ on $K$ such that $\pitil_1$ is indeed obtained by extending
$\pi_1$ to $\M(C_0(G_1))$. Put $x = (\io \ot \pi_1)(X)$. We have to
prove that $\cW := (\be \ot \io)(y)\;  x_{23}$ is a corepresentation
of $W$. Because $\pitil_1$ makes sense on all bounded Borel functions
on $\GGt$ (which is essentially the same thing as bounded Borel
functions on $G_1$), we get, by covariance $$\pitil_1(F(p_1(s \;
\cdot))) = y_s^* \pitil_1(F) y_s \; ,$$ for all bounded Borel
functions on $F$ on $G_1$. Integrating, we conclude that
$$\pitil_1((\om \ot \io)\al(F)) = (\om \ot \io)(y^* (1 \ot
\pitil_1(F))y) \; ,$$ for all bounded Borel
functions on $F$ on $G_1$ and $\om \in \B(H_2)_*$. From this, it
follows that Equation~\eqref{formula.3} holds for all $F \in C_0(G_1)$
and as we saw above, this yields that $\cW$ is a corepresentation of $W$.
\end{proof}

Observe that it follows immediately that the projection
$^*$-homomorphism of $S\uni$ onto $S$ corresponds exactly to the
projection $^*$-homomorphism of the full crossed product $G_2
\ltimes\full C_0(\GGt)$ onto the reduced crossed product $G_2
\ltimes\red C_0(\GGt)$. From the proof of the previous proposition, it
also follows that we have a non-degenerate $^*$-homomorphism $C_0(G_1)
\recht \M(S\uni)$.

Next, we consider covariant representations and we identify $S
\rtimes_{\text{\rm r,f}} \Sh$.

\begin{proposition} \label{bicroscovariant}
Suppose that $\cV$ is a representation and $\cW$ is a corepresentation
of $W$, such that $(\cV,\cW)$ is a covariant representation of $W$.

Take the unique corepresentations $x$ of $X$ and $y$ of $\Yh$
and the unique representations $a$ of $X$ and $b$ of $\Yh$ such that
$$\cV = b_{12} \; (\io \ot \al)(a) \quad\text{and}\quad \cW = (\be \ot
\io)(y) \; x_{23} \; .$$
Denote by $\pi_1$ and $\pi_2$ the representations of $C_0(G_1)$ and
$C_0(G_2)$ associated with $x$ and $b$, respectively.
\begin{enumerate}
\item The ranges
of $\pi_1$ and $\pi_2$ commute and we can define a non-degenerate
representation $\pi$ of $C_0(G_1 \times G_2)$. Extending first to $C_b(G_1
\times G_2)$ and then restricting to $C_0(G)$ using the map $(g,s)
\mapsto gs$, we
get a non-degenerate representation $\pitil$ of $C_0(G)$.
\item
The unitaries $a_g$ and $y_s$ commute for all $g \in G_1, s \in G_2$.
\item
The pair $(\pitil, (y_s \; a_g)_{(s,g)})$ is a covariant representation for the
action of $G_2 \times G_1$ on $C_0(G)$, where $G_2$ acts on the left and
$G_1$ on the right.
\item The images of the associated representations of $S
\rtimes\full \Sh$ and $(G_2 \times G_1) \ltimes\full C_0(G)$ coincide.
\end{enumerate}

Conversely, every covariant representation for the action of $G_2
\times G_1$ on $C_0(G)$ is obtained
in this way from a covariant pair for $W$.

In
particular, $$S \rtimes\full \Sh \cong (G_2 \times G_1) \ltimes\full
C_0(G)$$ in a natural way.
\end{proposition}
\begin{proof}
Let $(\cV,\cW)$ be a covariant pair for $W$ on $K$ and take, using
Proposition~\ref{universalcst}, $x,y,a$ and $b$ as in the statement of
the proposition. As we remarked right after
Definition~\ref{covariant}, both $\cW$ and $\cV$ are stably
isomorphic to $W$, but not jointly. This means that we have faithful,
normal $^*$-homomorphisms $\eta : M \recht \B(K)$ and $\etah : \Mh
\recht \B(K)$ such that $\cW = (\io \ot \eta)(W)$ and $\cV = (\etah
\ot \io)(W)$. Restricting to the von Neumann subalgebras
$\al(L^\infty(G_1))$ and $\cL(G_2) \ot 1$ of $M$ and
$\be(L^\infty(G_2))$ and $1 \ot \cR(G_1)$ of $\Mh$, we obtain
faithful, normal $^*\text{-}$homo\-morphisms $\pi_1,\pih_2,\pi_2$ and $\pih_1$
respectively, such that $$x = (\io \ot \pi_1)(X) \; , \quad y = (\io
\ot \pih_2)(\Yh) \; , \quad a = (\pih_1 \ot \io)(X)
\quad\text{and}\quad b = (\pi_2 \ot \io)(\Yh) \; .$$
The covariance of the pair $(\cV,\cW)$ is equivalent to each of the
following formulas: $$(\eta \ot \io)\de(z) = \cV(\eta(z) \ot 1) \cV^*
\quad\text{for all} \; z \in M \quad\quad (\io \ot \etah)\deh(z) =
\cW^* (1 \ot \etah(z)) \cW \quad\text{for all} \; z \in \Mh \; .$$
Evaluating these formulas on the generators of $M$ and $\Mh$ following
Equations~\eqref{eq.generatorM} and \eqref{eq.generatorMhat},
covariance is equivalent to each of the following two lines of formulas:
\begin{align}
& x_{12} \; (\io \ot \al)(X)_{134} = \cV_{234} \; x_{12} \;
\cV^*_{234} \quad , \quad y_{12} \; ((\io \ot \pi_1) \be \ot
\io)(\Yh)_{123} = \cV_{234} \; y_{12} \; \cV_{234}^* \quad , \label{formula.4}\\
& (\be \ot \io)(\Yh)_{124} \; b_{34} = \cW^*_{123} \; b_{34} \;
\cW_{123} \quad , \quad (\io \ot (\pi_2 \ot \io)\al)(X)_{234} \;
a_{34} = \cW^*_{123} \; a_{34} \; \cW_{123} \quad .\label{formula.5}
\end{align}
Using the explicit expression of $\cV$, the second formula of
Equation~\eqref{formula.4} becomes
\begin{equation} \label{formula.10}
b^*_{23} \; y_{12} \; ((\io \ot \pi_1) \be \ot
\io)(\Yh)_{123} \; b_{23} = (\io \ot \io \ot \al)(a_{23} \; y_{12} \;
a_{23}^*) \; .
\end{equation}
Because $\al(L^\infty(G_1)) \cap \cL(G_2) \ot 1 = \C$, we find a
unitary $v \in \B(H_2 \ot K)$, such that $a_{23} \; y_{12} \;
a_{23}^* = v_{12}$. Treating in the same way the second formula of
Equation~\eqref{formula.5}, we find that $y^*_{12} \; a_{23} \; y_{12}
\in 1 \ot \B(K \ot H_1)$. But, $y^*_{12} \; a_{23} \; y_{12} =
y^*_{12} \; v_{12} \; a_{23}$, which yields the existence of $u \in
\B(K)$ such that $v=y(1 \ot u)$. Because clearly both $v$ and $y$ are
corepresentations of $\Yh$, we find that $u=1$. Hence, we arrived at
the commutation of $a_{23}$ and $y_{12}$. This means that the
unitaries $a_g$ and $y_s$ commute for all $g \in G_1, s \in G_2$.

Next, we observe that $b^* (1 \ot F) b = (\pi_2 \ot \io)\de_2(F)$ for
all $F \in L^\infty(G_2)$. Combining this with the formula $(\de_2 \ot
\io)\al = (\io \ot \al)\al$, we can rewrite the first formula of
Equation~\eqref{formula.4} in the form
\begin{equation} \label{formula.11}
b_{23}^* \; x_{12} \; b_{23} = (\io \ot \io \ot \al)\bigl( a_{23} \;
x_{12} \; a_{23}^* \; (\io \ot (\pi_2 \ot \io)\al)(X^*) \bigr) \; .
\end{equation}
As above, it follows that $b_{23}^* \; x_{12} \; b_{23} \in \B(H_1 \ot
K) \ot 1$. Treating similarly the first formula of
Equation~\eqref{formula.5}, we find
\begin{equation} \label{formula.6}
x_{23} \; b_{34} \; x_{23}^* = (\be \ot \io \ot \io)\bigl(((\io \ot
\pi_1)\be \ot \io)(\Yh^*) \; y^*_{12} \;
b_{23} \; y_{12}  \bigr)
\end{equation}
and hence, that $x_{12} \; b_{23} \; x_{12}^*
\in 1 \ot \B(K \ot H_2)$. A similar reasoning as above yields now the
commutation of $x_{12}$ and $b_{23}$. This means that the ranges of
$\pi_1$ and $\pi_2$ commute. Combining this with
Equation~\eqref{formula.6}, we get
\begin{equation} \label{formula.7}
((\io \ot \pi_1)\be \ot \io)(\Yh) \; b_{23} = y^*_{12} \;
b_{23} \; y_{12} \; .
\end{equation}

Denote by $\pi$ the non-degenerate representation of $C_0(G_1 \times
G_2)$ on $K$ such that $\pi(F_1 \ot F_2)= \pi_1(F_1) \;
\pi_2(F_2)$. Dualizing the continuous map $(g,s) \mapsto gs$, we embed
$C_0(G)$ into $\M(C_0(G_1 \times G_2))$ and obtain the non-degenerate
representation $\pitil$ of $C_0(G)$ on $K$. We have to prove that
$(\pitil, (a_g \; y_s))$ is a covariant representation.

We claim that
$$\ga:=((\be \ot \io)\tau \ot \io)(\io \ot \de_2) : C_0(G_1 \times G_2)
\recht \M(\K(H_2 \ot H_1) \ot C_0(G_1 \times G_2))$$ is well-defined
and non-degenerate. It suffices to verify this statement on $C_0(G_1)
\ot 1$ and $1 \ot C_0(G_2)$ separately. For $C_0(G_1) \ot 1$, it
follows from Equation~\eqref{formula.non-deg}. On the other hand, we
observe that
$$(\ga \ot \io)(\Yh_{23}) = B_{12} \; B_{13} \; \Yh_{15} \; B^*_{13}
\; B^*_{12} \; \Yh_{45} \in \M(\K(H_2 \ot H_1) \ot C_0(G_1 \times G_2)
\ot \K(H_2)) \; ,$$
where $B \in \M(\K(H_2) \ot C_0(G_1))$ is again the canonical
implementation of $\be$. This proves our claim. Using this last
equation and Equation~\eqref{formula.7}, we also observe that
\begin{align*}
(\io \ot \io \ot \pi \ot \io)(\ga \ot \io)(1 \ot \Yh) & = B_{12} \;
((\io \ot \pi_1)\be \ot \io)(\Yh)_{134} \; B_{12}^* \; b_{34} \\
& = (\be \ot \io \ot \io)(y^*_{12} \; b_{23} \; y_{12}) = (\be \ot
\io)(y^*)_{123} \; (\pi \ot \io)(1 \ot \Yh)_{34} \; (\be \ot
\io)(y)_{123} \; .
\end{align*}
Combining this with Equation~\eqref{formula.3}, which holds because
$\cW$ is a corepresentation of $W$, we get
\begin{equation} \label{formula.9}
(\io \ot \io
\ot \pi)\ga(F) = (\be \ot \io)(y^*) \; (1 \ot 1 \ot \pi(F)) \; 
(\be \ot \io)(y) \; ,
\end{equation}
for all $F \in C_0(G_1 \times G_2)$. By
non-degenerateness, the same holds for $F \in C_0(G)$, which exactly
means that
$$\pitil(F(s \; \cdot)) = y_s^* \; \pitil(F) \; y_s \; ,$$
for all $F \in C_0(G)$ and $s \in G_2$.

One similarly proves that $\pitil(F(\cdot\; g)) = a_g \; \pitil(F) \;
a_g^*$ for all $F \in C_0(G)$ and $g \in G_1$. So,
$(\pitil, (a_g \; y_s))$ is a covariant representation.

Denote the image of the representation of $S \rtimes\full \Sh$ associated with
$(\cV,\cW)$ by $A$. Then, by definition
$$A=[(\io \ot \om)(\cV) \; (\mu \ot \io)(\cW) \mid \om,\mu \in \B(H_2
\ot H_1)_* ] \; .$$
It follows from Proposition~\ref{universalcst} and its proof that
$$A=[\pih_2(C^*(G_2) \; \pi_1(C_0(G_1)) \; \pi_2(C_0(G_2)) \; 
\pih_1(C^*(G_1)) ] = [\pih_2(C^*(G_2) \; \pi(C_0(G_1 \times G_2)) \; 
\pih_1(C^*(G_1)) ] \; .$$
From Equation~\eqref{formula.9}, we know that
$$\pi\bigl( (\om \ot \io \ot \io)(\tau \ot \io)(\io \ot \de_2)(F)
\bigr) = (\om \ot \io)(y^* \;  (1 \ot \pi(F)) \; y)$$
for all $F \in C_0(G_1 \times G_2)$ and $\om \in \B(H_2)_*$. 
Lemma~\ref{nogeen} following this proof, tells us that the left hand side
spans a dense subset of $\pitil(C_0(G))$. So, we get
$$[\pih_2(C^*(G_2) \; \pitil(C_0(G)) \; \pih_1(C^*(G_1) ] = 
[\pih_2(C^*(G_2) \; \pi(C_0(G_1 \times G_2)) \; \pih_1(C^*(G_1) ] \;
.$$
Because both $A$ and the image of $(G_2 \times G_1) \ltimes\full C_0(G)$
are \cst-algebras, we easily get that $A$ indeed equals the
image of $(G_2 \times G_1) \ltimes\full C_0(G)$.

Suppose now that, conversely, we have a covariant representation
$(\pitil, (y_s \; a_g))$ for the action of $G_2 \times G_1$ on $C_0(G)$.
The representation $\pitil$ gives rise to a measure class on $G$ which
is invariant by multiplication on the left by $G_2$ and on the right
by $G_1$. We claim that this measure class is supported by $G_1G_2$.
Because the transformation $(s,x) \mapsto (s,sx)$ on
$G_2 \times G$ preserves Borel sets of measure zero, it is enough to
prove that for any $x \in G$, the set of all $s \in G_2$ such that $sx
\not\in G_1G_2$ has measure zero. But this is the case, as we already
saw in the proof of Proposition~\ref{universalcst}.

Hence, the Borel functional calculus provides a non-degenerate representation $\pi$
of $C_0(G_1 \times G_2)$ such that $\pitil$ is obtained by extending
$\pi$ to $\M(C_0(G_1 \times G_2))$ and then restricting to $C_0(G)$
through the map $(g,s) \mapsto gs$. We
denote by $\pi_1$ and $\pi_2$ the corresponding representations of
$C_0(G_1)$ and $C_0(G_2)$. Then, we can define $x=(\io \ot \pi_1)(X)$
and $b = (\pi_2 \ot \io)(\Yh)$. Defining
$$\cV = b_{12} \; (\io \ot \al)(a) \quad\text{and}\quad \cW = (\be \ot
\io)(y) \; x_{23} \; ,$$
it follows from Proposition~\ref{universalcst} that $\cV$ is a
representation of $W$ and $\cW$ is a corepresentation of $W$.

We know that $\pitil(F(s \; \cdot)) = y_s^* \pitil(F) y_s$ for all $s
\in G_2$ and all bounded Borel functions $F$ on $G$. Integrating, it
follows in particular that
$$\pitil\bigl( (\om \ot \io \ot \io)(\be \ot \io)\de_2(F) \bigr) =
(\om \ot \io)(y^*(1 \ot \pi_2(F))y)$$
for all $\om \in \B(H_2)_*$ and all $F \in C_0(G_2)$. Evaluating this formula on $F=(\io \ot \mu)(\Yh)$, we get
$$((\io \ot \pi_1)\be \ot \io)(\Yh) \; b_{23} = y_{12}^* \;
b_{23} \; y_{12} \; ,$$
which makes sense because $\be : C_0(G_2) \recht \M(\K(H_2) \ot
C_0(G_1))$ is well-defined and non-degenerate (using the canonical
implementation of $\be$). Because $a_{23}$ and
$y_{12}$ commute, it follows that Equation~\eqref{formula.10} holds
and hence, also the second formula of Equation~\eqref{formula.4}
holds.

Because $\pitil(F(\cdot\;g)) = a_g \pitil(F) a_g^*$ for all $g \in
G_1$ and all bounded Borel functions $F$ on $G$, we find similarly
that
$$\pitil\bigl( (\io \ot \io \ot \om)(\io \ot \al)\de_1(F) \bigr) =
(\io \ot \om)(a(\pi_1(F) \ot 1) a^*)$$
for all $\om \in \B(H_1)_*$ and all $F \in C_0(G_1)$, which yields, in
the same way,
$$x_{12} \; (\io \ot (\pi_2 \ot \io)\al)(X) = a_{23} \; x_{12} \;
a_{23}^* \; ,$$
where $(\pi_2 \ot \io)\al$ makes again sense using the canonical
implementation of $\al$.
Because $x_{12}$ and $b_{23}$ commute, we get that
Equation~\eqref{formula.11} holds and hence, also the first formula in
Equation~\eqref{formula.4}. Combining both formulas of
Equation~\eqref{formula.4} and the definitions of $\cV$, $\cW$ and $W$, we get
\begin{align*}
\cV_{345} \; \cW_{123} \; \cV_{345}^* &= (\be \ot \io)(y)_{123} \;
((\io \ot \io \ot \pi_1)(\io \ot \de_1)\be \ot \io)(\Yh)_{1234} \;
x_{23} \; (\io \ot \al)(X)_{245} \\
&= \cW_{123} \; (\be \ot \io)(\Yh)_{124} \; (\io \ot \al)(X)_{245} = 
\cW_{123} \; W_{1245} \; ,
\end{align*}
where we used the formula $(\io \ot \pi_1)\de_1(F) = x (F \ot 1) x^*$
for all $F \in C_0(G_1)$. So, we precisely get that $(\cV,\cW)$ is a
covariant pair for $W$.
\end{proof}

The following lemma was needed to prove the previous proposition. It
is completely analogous to Lemma~\ref{quotient} above.

\begin{lemma} \label{nogeen}
Using the map $(g,s) \mapsto gs$, we embed $C_0(G) \recht \M(C_0(G_1
\times G_2)) \subset L^\infty(G_1 \times G_2)$. Then, we have
$$C_0(G) = [(\om \ot \io \ot \io)(\tau \ot \io)(\io \ot \de_2)(C_0(G_1
\times G_2)) \mid \om \in \B(H_2)_* ]  \; .$$
\end{lemma}
\begin{proof}
We have to prove that the functions
$$gt \mapsto \int K_2(s) F_1(p_1(sg)) F_2(p_2(sg)t) \; ds$$
belongs to $C_c(G)$ whenever $F_1 \in C_c(G_1), K_2,F_2 \in C_c(G_2)$
and that they span a dense subset of $C_0(G)$. It suffices to consider
$F_1$ of the form $\tilde{K}_1 * P_1$ with $K_1,P_1 \in C_c(G_1)$ and
$$F_1(g) = \int K_1(h) P_1(hg) \; dh \; .$$
Then, we can define bounded Borel functions $K$ and $P$ with compact
support on $G$ by the formulas
$$P(gs) = P_1(g) F_2(s) \quad , \quad K(gs) = K_1(g) K_2(s)
\sde_1(g)^{-1} \sde(g)$$
and $P$ and $K$ equal $0$ outside $G_1G_2$.
Because $K$ and $P$ belong to $L^2(G)$ and have compact support,
the function $H:= K * P$ defined by
$$H(x) = \int K(y) P(yx) \; dy$$ belongs to $C_c(G)$. But, using Proposition~\ref{density},
\begin{align*}
H(gt) &= \int K(y) P(ygt) \; dy = \iint K(hs) P(hsgt) \sde_1(h)
\sde(h)^{-1} \; dh \; ds \\
&= \iint K_1(h) K_2(s) P_1(hp_1(sg)) F_2(p_2(sg)t) \; dh \; ds = \int
K_2(s) F_1(p_1(sg)) F_2(p_2(sg)t) \; ds \; .
\end{align*}
This ends the proof, because it is clear that the functions $K$ and
$P$ span dense subspaces of $L^2(G)$ and hence the functions $H$ span
a dense subspace of $C_0(G)$.
\end{proof}

We also characterize the reduced crossed products in the following easily proved
proposition.
\begin{proposition} \label{reducedcovariant}
The regular covariant pair $(W,W)$ for $W$ corresponds through the
procedure of Proposition~\ref{bicroscovariant} to a
representation of $(G_2 \times G_1)
\ltimes\full C_0(G)$ which is stably isomorphic to the regular representation. In particular,
$$S \rtimes\red \Sh \cong (G_2 \times G_1) \ltimes\red C_0(G)$$
and the projection $^*$-homomorphisms from the full onto the reduced
crossed products are intertwined by the isomorphisms $S
\rtimes_{\text{\rm f,r}} \Sh \cong (G_2 \times G_1) \ltimes_{\text{\rm f,r}} C_0(G)$.
\end{proposition}
\begin{proof}
Because the orbit of $e$ under the action of $G_2 \times G_1$, which
is $G_2 G_1$, is dense in $G$, it follows that the covariant
representation associated with this orbit, on the Hilbert space
$L^2(G_2 \times G_1)$, is stably isomorphic to the regular
representation of $(G_2 \times G_1)
\ltimes\full C_0(G)$. It is clear that the $\cV$ and $\cW$ corresponding
to this covariant representation are twice $W$.
\end{proof}

Using this proposition, we prove our main theorem.

\begin{theorem} \label{reg-and-semireg}
The multiplicative unitary $W$ of the bicrossed product l.c.\ quantum
group $(M,\de)$ is regular if and only if the map
$$\te : G_1 \times G_2 \recht G : \te(g,s) = gs$$
is a homeomorphism of $G_1 \times G_2$ onto $G$. The multiplicative
unitary $W$ is semi-regular if and only if $\te$ is a homeomorphism of
$G_1 \times G_2$ onto an open subset of $G$ with complement of measure
zero.
\end{theorem}
\begin{proof}
In the proof of the previous proposition, we have seen that $S
\rtimes\red \Sh$ is precisely given by the image, denoted by $A$, of the irreducible representation of $(G_2 \times
G_1) \ltimes\full C_0(G)$ corresponding to the free orbit
$G_2G_1$. From Lemma~\ref{hulplemma}, it follows that $A = \K(H_2 \ot H_1)$ if and only if this orbit
is closed and homeomorphic to $G_2 \times G_1$ and that $\K(H_2 \ot
H_1) \subset A$ if and only if this orbit is locally closed and
homeomorphic to $G_2 \times G_1$. Because the orbit is dense and
because of Proposition~\ref{char}, we
precisely arrive at the statement of the theorem, using also that $G_1
G_2 = (G_2 G_1)^{-1}$.
\end{proof}

In the next section, we will give examples were the image of $\te$ is
not open and hence, the associated multiplicative unitary is not semi-regular.
Nevertheless, it should be observed that $W$, being the regular
representation of a l.c.\ quantum group, is always manageable in the
sense of \cite{Wor}. So, not all manageable multiplicative unitaries
are semi-regular.

The next lemma is well known (see \cite{fack} for a related result), but
we include a short proof for completeness.

\begin{lemma} \label{hulplemma}
Suppose that a l.c.\ group $G$ acts continuously (on the right) on a l.c.\ space
$X$. Let $x_0 \in X$ have a free orbit, i.e.\ the map $\te: G \recht X
: \te(g) = x_0 \cdot g$ is injective.

Denote by $\pi$ the representation of $C_0(X) \rtimes\full G$ on
$L^2(G)$ corresponding to the orbit of $x_0$. Then,
\begin{enumerate}
\item The image of $\pi$ contains $\K(L^2(G))$ if and only if the
  orbit $\te(G)$ is locally closed and $\te$ is a homeomorphism of $G$
  onto $\te(G)$.
\item The image of $\pi$ is equal to $\K(L^2(G))$ if and only if the
  orbit $\te(G)$ is closed and $\te$ is a homeomorphism of $G$ onto
  $\te(G)$.
\end{enumerate}
\end{lemma}
\begin{proof}
Suppose first that the image of $\pi$ contains $\K(L^2(G))$. Equip
$\te(G)$ with its relative topology and assume that $\te^{-1}$ is not
continuous. Then, we find elements $g_n \in G$ such that $g_n$ $(n
\geq 1)$ remains outside a neighborhood of
$g_0$, but $x_0 \cdot g_n \recht x_0 \cdot g_0$. Consider
the dense subalgebra of the image of $\pi$ consisting of the operators
$\gamma(F)$, $F \in C_c(G \times X)$, defined by
$$\bigl( \gamma(F) \xi)(g) = \int_G F(h,x_0 \cdot g) \; \xi(gh) \; dh
\; .$$
Take a function $\eta \in L^2(G)$ with $\|\eta\|=1$ and with small
enough support such that $\langle \lambda_{g_n} \eta, \lambda_{g_0}
\eta \rangle = 0$ for all $n \geq 1$, where $(\lambda_g)$ is the left
regular representation of $G$. One verifies immediately that
$$\langle \lambda_{g_n} \eta , \gamma(F) \lambda_{g_n} \eta \rangle
\recht \langle \lambda_{g_0} \eta , \gamma(F) \lambda_{g_0} \eta
\rangle \quad\text{for all}\quad F \in C_c(G \times X) \; .$$
Hence, the same holds for all $a \in \pi(C_0(X) \rtimes\full G)$
instead of $\gamma(F)$ and, in particular, for all $a \in
\K(L^2(G))$. This gives a contradiction when we take for $a$ the
projection on $\lambda_{g_0} \eta$. So, $\te^{-1}$ is continuous. This
means that $\te(G)$ is locally compact in its relative topology, which
precisely means that $\te(G)$ is locally closed in $X$.

Suppose next that the image of $\pi$ is precisely $\K(L^2(G))$. Denote
by $Y$ the closure of $\te(G)$. From the previous paragraph, we
already know that $\te(G)$ is open in $Y$ and that $\te$ is a
homeomorphism. Suppose that $\te(G) \neq Y$, take $x_1 \in Y
\setminus \te(G)$ and take $g_n \in G$ such that $x_0\cdot g_n \recht
x_1$. Then, $g_n$ goes to infinity. Writing $\pi_1$ for the representation of $C_0(X) \rtimes\full G$ on
$L^2(G)$ corresponding to the orbit of $x_1$, whose image contains the
dense subalgebra of operators $\gamma_1(F)$, $F \in C_c(G \times
X)$, we observe that
$$\langle \lambda_{g_n} \eta, \gamma(F) \lambda_{g_n} \eta \rangle
\recht \langle \eta, \gamma_1(F) \eta \rangle \quad\text{for all}\quad
\eta \in L^2(G), F \in C_c(G \times X) \; .$$
But, $\lambda_{g_n} \eta \recht 0$ weakly, as $g_n \recht \infty$. Because the
image of $\pi$ is supposed to be exactly the compact operators, it
follows from the previous formula that $\langle \eta, \gamma_1(F)
\eta \rangle = 0$ for all $F$ and $\eta$. This is a contradiction.

The converse implications are easy to prove.
\end{proof}

Using the theory that we developed so far, it is also easy to give examples
of l.c.\ quantum groups such that the projection $^*$-homomorphism
from $S \rtimes\full \Sh$ onto $S \rtimes\red \Sh$ is not
faithful. Loosely speaking, this means that the action of $(S,\de)$ on
itself by translation is not amenable and in particular, not proper (whatever this means).

\begin{proposition} \label{conjugated}
Suppose that $G_1$ and $G_2$ are conjugated, i.e.\ there exists an
element $z_0 \in G$ such that $G_2 = z_0 G_1 z_0^{-1}$.
\begin{enumerate}
\item $S \rtimes_{\text{\rm f,r}} \Sh$ is strongly Morita equivalent
  with $S_{\text{\rm u,r}}$.
\item The projection $^*$-homomorphism
from $S \rtimes\full \Sh$ onto $S \rtimes\red \Sh$ is faithful if and
only if $G_1$ is amenable.
\item $M \cong \Mh \cong \B(L^2(G_1))$.
\end{enumerate}
\end{proposition}
\begin{proof}
{\it a)} We use an argument which is essentially contained in Rieffel's
paper \cite{Rie}, but we include a sketch of it for completeness. Write $B=C_0(G/G_1)$. On $C_c(G)$, the continuous compactly
supported functions on $G$, we define a $B$-valued inner product by
$$\langle \xi ,\eta \rangle(x) = \int_{G_1} (\overline{\xi} \eta)(x
g^{-1}) \; dg \; .$$
Completion yields the Hilbert $B$-module $\cE$. Because the right
action of $G_1$ on $G$ is proper, there is only one crossed product
$C_0(G) \rtimes G_1$, which can be identified with $\K(\cE)$, the
\lq compact\rq\ operators on the Hilbert $B$-module $\cE$. Because
$\cE$ is full, we have a Morita equivalence between $C_0(G) \rtimes
G_1$ and $B$. Then, we have
$$S \rtimes_{\text{\rm f,r}} \Sh \cong (G_2 \times G_1) \ltimes_{\text{\rm f,r}} C_0(G) \cong G_2
\ltimes_{\text{\rm f,r}} \K(\cE) \cong \K(G_2 \ltimes_{\text{\rm f,r}}
\cE) \; , $$
where all isomorphisms are natural and intertwine the projection of
full onto reduced crossed products and where $G_2 \ltimes_{\text{\rm
    f,r}} \cE$ is the obvious full Hilbert $G_2 \ltimes_{\text{\rm
    f,r}} C_0(G/G_1)$-module. In the identifications above, only the
isomorphism $G_2 \ltimes\full \K(\cE) \cong \K(G_2 \ltimes\full \cE)$
requires some care: consider the \cst-algebra $\K(\cE \oplus B)$ in
which $\K(\cE)$ is a full corner and on which $G_2$ acts by
automorphisms. Then, $G_2 \ltimes\full \K(\cE)$ is a full corner of
$G_2 \ltimes\full \K(\cE \oplus B)$ and is as such identified with
$\K(G_2 \ltimes\full \cE)$.

Observing that the right multiplication by
$z_0^{-1}$ gives a homeomorphism of $G/G_1$ onto $G/G_2$ intertwining
the left action of $G_2$, we see that $G_2 \ltimes_{\text{\rm
    f,r}} C_0(G/G_1) \cong S_{\text{\rm u,r}}$. So, we have proven the
required strong Morita equivalence.

{\it b)}
Because of the Morita equivalence above, the projection of $S
\rtimes\full \Sh$ onto $S \rtimes\red \Sh$ is faithful if and only if
the projection of $S\uni$ onto $S$ is faithful. From Theorem~15 in
\cite{DQV}, it follows that this last projection
is faithful if and only if $G_2$ is amenable. Because
of conjugacy, this is equivalent to the amenability of $G_1$.

{\it c)} As above, we observe that the left action of $G_2$ on $G/G_2$
is isomorphic to the left action of $G_2$ on $G/G_1$. Using the
isomorphism $L^\infty(G/G_1) \cong L^\infty(G_2)$, we see that the
action of $G_2$ on $L^\infty(G/G_2)$ is isomorphic to the left
translation of $G_2$ on $L^\infty(G_2)$. Hence, $M \cong \B(L^2(G_2))$.
\end{proof}

\begin{remark}
Although in the situation of Proposition~\ref{conjugated}, we have a strong Morita equivalence between $S
\rtimes\full \Sh$ and $S\uni$ and hence, an isomorphism $S
\rtimes\full \Sh \ot \K \cong S\uni \ot \K$, this isomorphism is very
much \lq twisted\rq\ for the following reason: if $G_1$ is
non-amenable, $S\uni$ is very different from $S$, but nevertheless, we
claim that, for any l.c.\ quantum group $(M,\de)$, we have a natural, injective
$^*$-homomorphism $M \recht \M(S \rtimes\full \Sh)$. As usual we write
$S$ and $\Sh$ for the underlying \cst-algebras as in Equation~\eqref{eq.SandSh}.

From the remark after Definition~\ref{covariant}, it follows that we
can realize (in a natural way) $S \rtimes\full \Sh$ on a Hilbert space
$K$ such that there exist normal, faithful $^*$-homomorphisms $\pi :
M \recht \B(K)$ and $\pih : \Mh \recht
\B(K)$ such that $S \rtimes\full \Sh = [\pi(S) \; \pih(\Sh)]$. We
prove that $\pi(M) \subset \M(S \rtimes\full \Sh)$. Take
$x \in S^{\prime\prime}$ and observe that
$$(\pih \ot \io)(V) (\pi(x) \ot 1) = (\pi \ot \io)\de(x) (\pih \ot \io)(V) =
(\pi \ot \io)(W^*) (1 \ot x) (\pi \ot \io)(W) (\pih \ot \io)(V) \; ,$$
where $V \in \M(\Sh \ot \K)$ is the right regular representation of
$(M,\de)$ and $W \in \M(S \ot \K)$ is the left regular
representation. Hence,
\begin{align*}
(S \rtimes\full \Sh) \pi(x) &= [\pi(S) (\io \ot \om) \bigl( 
(\pi \ot \io)(W^*) (1 \ot x) (\pi \ot \io)(W) (\pih \ot \io)(V) \bigr)
\mid \om \in \B(H)_* ] \\ &= [(\io \ot \om) \bigl( (\pi(S) \ot x) (\pi \ot \io)(W) (\pih \ot \io)(V) \bigr)
\mid \om \in \B(H)_* ] \subset S \rtimes\full \Sh \; .
\end{align*}
So, we are done.
\end{remark}

\section{Examples} \label{sec.examples}

We start off by presenting a fairly general example. Suppose that
$\cA$ is a l.c.\ ring, with unit, but not necessarily commutative and $\cA \neq \{0\}$. Denote by $\cAu$ the group of invertible
elements in $\cA$. Then, $\cAu$ is a l.c.\ group by considering it as
the closed subspace $\{ (a,b) \mid ab=1=ba \}$ of $\cA \times \cA$.

Next, we can define the $ax+b$-group on this l.c.\ ring $\cA$.
$$G= \cAu \times \cA \quad\text{with}\quad (a,x)(b,y) = (ab,x+ay) \;
.$$
Defining the closed subgroups (isomorphic with $\cAu$)
$$G_1 = \{ (a,a-1) \mid a \in \cAu \} \quad\text{and}\quad 
G_2 = \{ (b,0) \mid b \in \cAu \} \; ,$$ we observe that $G_1 G_2$
consists of all pairs $(a,x)$ such that $x+1 \in \cAu$. Hence, we get
a matched pair of l.c.\ groups if and only if $\cAu$ has complement of
(additive) Haar measure zero in $\cA$. Observe that $G_1$ and $G_2$
are conjugated by the element $(-1,-1) \in G$.

Suppose that $\cA$ is a l.c.\ ring,
different from $\{ 0 \}$, with $\cA \setminus \cA^*$ of (additive) Haar measure zero.

Denote the bicrossed product multiplicative unitary associated to the
above matched pair by
$\cW_\cA$ and the associated l.c.\ quantum group by $(M,\de)_\cA$. 

\begin{proposition} 
We have the following properties.
\begin{itemize}
\item
The multiplicative unitary $\cW_\cA$ is not regular. It is
semi-regular if and only if $\cAu$ is open in $\cA$.
\item
The projection of $S \rtimes\full \Sh$ onto $S \rtimes\red \Sh$ is
faithful if and only if the l.c.\ group $\cAu$ is amenable. 
\item
The dual $(\Mh,\deh)_\cA$ is isomorphic to the opposite quantum group $(M,\de\Op)_\cA$.
\item
The associated \cst- and von Neumann algebras are given by
\begin{align*}
& S_{\text{\rm u,r}} \cong \Sh_{\text{\rm u,r}} \cong \cAu \ltimes_{\text{\rm f,r}} C_0(\cA)
\quad\text{where $\cAu$ multiplies $\cA$ on the left} \; , \\
& M \cong \Mh \cong \B(L^2(\cAu)) \; ,\\
& S \rtimes_{\text{\rm f,r}} \Sh \cong \K(L^2(\cAu)) \ot \bigl(\cAu
\ltimes_{\text{\rm f,r}} C_0(\cA)\bigr) \; .
\end{align*}
\end{itemize}
\end{proposition}
\begin{proof}
Using Theorem~\ref{reg-and-semireg} and Proposition~\ref{conjugated},
we get immediately the first two statements of the proposition.

To prove the third statement, define the unitary $\cU$ on $L^2(\cAu
\times \cAu)$ by the formula $$(\cU \xi)(b,a) = \sde_{\cAu}(ab)^{1/2}
\xi(a^{-1},b^{-1})  \; ,$$ where $\sde_{\cAu}$ is the modular function of
$\cAu$. A straightforward calculation on the generators shows that $\cU M \cU^* = \Mh$
and that this isomorphism intertwines $\de\Op$ on $M$ and $\deh$ on $\Mh$

There is a homeomorphism from $G/G_2$ onto $\cA$ mapping
$\overline{(a,x)}$ to $x$, intertwining the left action of $G_2$ on
$G/G_2$ with the left multiplication of $\cAu$ on $\cA$.  So, 
using Propositions~\ref{cst} and \ref{universalcst}, we get $S_{\text{\rm u,r}}
\cong \cAu \ltimes_{\text{\rm f,r}} C_0(\cA)$ where $\cAu$ multiplies
$\cA$ on the left. Because of the already proven third statement, we
have $S_{\text{\rm u,r}} \cong \Sh_{\text{\rm u,r}}$ and so, we have
proven the first formula of the fourth statement.

The second formula of the fourth statement follows from Proposition~\ref{conjugated}.

Finally, using the homeomorphism $(d,x) \mapsto (d,x-d)$ of $G$, the
action of $G_2 \times G_1$ on $G$ is isomorphic to the action $b \cdot
(d,x) \cdot a = (bda,bx)$ of $G_2 \times G_1$ on $G$. If we take the
crossed product with the action of $a \in \cAu$, we obtain
$\K(L^2(\cAu)) \ot C_0(\cA)$, on which $\cAu$ acts diagonally. This
diagonal action is isomorphic with the amplified action on $C_0(\cA)$
and we obtain the third formula, using Propositions~\ref{bicroscovariant} and \ref{reducedcovariant}.
\end{proof}

Observe that the von Neumann algebraic picture of $(M,\de)_\cA$ is
fairly trivial, because $M \cong \B(L^2(\cAu))$. Nevertheless, the
\cst-algebra picture is far from trivial. As we will see in concrete
examples below, the \cst-algebra $S$ need not be type I.

An example of a non-amenable $\cA^*$ is, for instance, given by $\cA
= M_2(\R)$. In that case the projection of $S \rtimes\full \Sh$ onto
$S \rtimes\red \Sh$ is not an isomorphism.

\begin{example} \label{ex.qp}
Let $\cP$ be an infinite set of prime numbers such that $$\sum_{p \in \cP}
\frac{1}{p} < \infty \; .$$ Define $$\cA = {\prod_{p \in
  \cP}}' \; \Q_p \; ,$$ where $\Q_p$ is the l.c.\ field of $p$-adic numbers
and where the prime means that we take the restricted Cartesian
product: we only consider elements $(x_p)$
that eventually belong to the compact sub-ring $\Z_p$ of $p$-adic
integers. When we equip $\cA$ with the usual l.c.\ topology, we
obtain a l.c.\ ring such that $\cAu$ has complement of measure zero in
$\cA$ but $\cAu$ has empty interior in $\cA$. The fact that $\cAu$ has
complement of measure zero follows from a straightforward application
of the Borel-Cantelli lemma: normalizing the Haar measure on $\Q_p$
such that $\Z_p$ has measure one, we observe that $\Z_p \setminus \Z_p\unit$
has measure $\frac{1}{p}$, which is assumed to be summable over $p \in \cP$.

Observe that $\cA$ is a (restricted) ring of adeles and that $S_\cA \cong
\cAu \ltimes C_0(\cA)$. Such a \cst-algebra $S_\cA$ is not of type I and
is intensively studied by J.B.~Bost and A.~Connes in
\cite{BostC}\footnote{More precisely, they study the Hecke algebra,
  which is $p (\cAu \ltimes C_0(\cA)) p$ for a natural projection $p$.}. 
\end{example}

\begin{remark}
Observe that in the previous example, we can replace $\Q_p$ by the
field of formal Laurent series over the finite field $F_p$ with $p$
elements. We have the compact sub-ring of formal power series over
$F_p$ and can perform an analogous construction.
\end{remark}

\begin{example}
Another, more canonical example can be given, much in the same spirit
as Example~\ref{ex.qp}. Denote for every prime number $p$ by $K_p$ the
(unique) non-ramified extension of degree 2 of $\Q_p$. Let $\cO(K_p)$
be the compact ring of integers in $K_p$. Because the extension is
non-ramified, the Haar measure of $\cO(K_p) \setminus \cO(K_p)\unit$
is $\frac{1}{p^2}$. Because this sequence is summable, we can take
$$\cA = {\prod_{p \; \text{prime}}}\hspace{-.17cm}\raisebox{0.15cm}[0cm][0cm]{$^\prime$}\hspace{.1cm} K_p \; .$$
\end{example}

\begin{remark} \label{strange.mult.un}
Consider the easiest case $\cA = \R$. Define $K = L^2(\R^*)$. A
covariant representation on the Hilbert space $K$ for $G_2 \times G_1$
acting on $C_0(G)$ is given by
$$\pi(F)(x) = F(x,x) \quad\text{for} \; F \in C_0(G) ,x \in \R^*
\qquad\text{and}\qquad (\lambda_b \rho_a)_{(b,a) \in G_2 \times G_1}
\; .$$
We can then take the covariant image of the bicrossed product
multiplicative unitary $W_\R$ and this yields the multiplicative
unitary $\cW$ on $K$ given by
$$(\cW \xi)(x,y) = \xi(\frac{xy}{x+y+1}, \frac{y}{x+1}) \; .$$
Observing that $(\cW^* \xi)(x,y) = \xi(\frac{x(y+1)}{y},x+y+xy)$, we
see that both $\cW$ and $\cW^*$ leave $L^2(\R^*_+) \ot
L^2(\R^*_+)$ invariant. So, the restriction of $\cW$ gives us a multiplicative
unitary $\tilde{\cW}$ on $K_+:=L^2(\R^*_+)$. It is easy to check that the
weak closure of $S_{\tilde{\cW}}$ consists precisely of all the
operators $T \in \B(K_+)$ such that $P_{[x,+\infty[} K_+$ is an
invariant subspace of $T$ for all $x \in \R^*_+$. Hence, this weak closure is
not invariant under involution. This means that $S_{\tilde{\cW}}$ is
not a \cst-algebra. Analogously, $\Sh_{\tilde{\cW}}$ is not a
\cst-algebra.

The multiplicative unitary $\tilde{\cW}$ is very singular and should
certainly not be considered as the multiplicative unitary of a quantum
group.

Using the transformation $v: ]0,1[ \recht R^*_+ : v(x) =
\frac{x}{1-x}$, the multiplicative unitary $\tilde{\cW}$ is
transformed to a multiplicative unitary $X$ on $L^2(]0,1[,d\mu)$ given
by
$$(X \xi)(x,y) = \xi(xy,y \frac{1-x}{1-xy}) \; .$$
The transformation $v(x,y) = (xy, y \frac{1-x}{1-xy})$ is studied in
\cite{kashaev} and is shown to be essentially the only pentagonal
transformation on $]0,1[$ of the form $v(x,y) = (xy,u(x,y))$ with $u$
continuously differentiable.

We can construct an even more singular multiplicative unitary $Y$ on the
Hilbert space $l^2(\Q^*_+)$ which is formally given by the same
formula as $\tilde{\cW}$. It is easy to check that $S_Y = [
\theta_{e_s} \theta_{e_r}^* \mid 0 < r < s, r,s \in \Q ]$, where
$(e_q)$ denotes the obvious basis in $l^2(\Q^*_+)$. Instead of $\Q$,
we can as well take any countable subfield of $\R$. 
\end{remark}

\begin{example}
Up to now, we considered in fact only one example: the $ax+b$ group
over an arbitrary l.c.\ ring with two natural subgroups. If $\cA =
\R$, this example is the easiest non-trivial example of a matched pair
of real Lie groups (and the only non-trivial example when $G$ is of
dimension 2).
But, there are a lot of examples of matched pairs of algebraic groups, even
in low dimensions, see e.g.\ \cite{VV2}. Often, the field $\R$ or $\C$
can be replaced by any l.c.\ ring $\cA$ with $\cA \setminus \cAu$ of
measure zero and $\cA \neq \{0\}$. Fix such a l.c.\ ring $\cA$.

We give two examples. First, define $G =
GL_2(\cA)$, the two by two matrices over $\cA$. Take
$$G_1 = \begin{pmatrix} \cAu & \cA \\ 0 & 1 \end{pmatrix}
\quad\text{and}\quad G_2 = \begin{pmatrix} 1 & 0 \\ \cA & \cAu
\end{pmatrix} \; .$$
Observe that $G_1$ and $G_2$ are conjugated by
$\bigl(\begin{smallmatrix} 0 & 1 \\ 1 & 0 \end{smallmatrix}
\bigr)$.

Next, we give an example in the same spirit as \cite{VV}, Section
5.4. Fix $q \in \cA$ and suppose $q$ to be central. Define on
$\cB_q:=\cA^4$ the structure of a l.c.\ ring by putting
$$\begin{pmatrix} a & b \\ c & d \end{pmatrix}\subq + 
\begin{pmatrix} a' & b' \\ c' & d' \end{pmatrix}\subq =
\begin{pmatrix} a+a' & b+b' \\ c+c' & d+d' \end{pmatrix}\subq
\qquad\text{and}\qquad \begin{pmatrix} a & b \\ c & d \end{pmatrix}\subq \cdot
\begin{pmatrix} a' & b' \\ c' & d' \end{pmatrix}\subq =
\begin{pmatrix} aa' + qb c' & ab'+bd' \\ ca'+dc' & dd' + qcb'
\end{pmatrix}\subq \; .$$
Defining $\pi_q \left(\begin{smallmatrix} a & b \\ c & d
  \end{smallmatrix}\right)\subqq = \left(\begin{smallmatrix} a & b \\ qc & d
  \end{smallmatrix}\right) \oplus \left(\begin{smallmatrix} a & qb \\ c & d
  \end{smallmatrix}\right)$, we identify
$\cB_q$ with the closed subring of
$M_2(\cA) \oplus M_2(\cA)$ consisting of the elements
$\left(\begin{smallmatrix} a & b \\ c & d \end{smallmatrix}\right)
\oplus \left(\begin{smallmatrix} a & b' \\ c' & d
  \end{smallmatrix}\right)$ satisfying $b'=q b$ and $c=qc'$. If $q \in
\cAu$, the first component of $\pi_q$ gives an isomorphism $\cB_q
\cong M_2(\cA)$.

Define $G^q = (\cB_q)\unit$ and
$$G_1^q = \{ \begin{pmatrix} a & b \\ 0 & d \end{pmatrix}\subq \mid a,d
\in \cAu, b \in \cA \} \qquad\text{and}\qquad G_2^q = \{
\begin{pmatrix} 1 & 0 \\ x & 1 \end{pmatrix}\subq \mid x \in \cA \} \;
.$$
Observe that $G_1^q \cong G_1^1 \cong$ the group of upper triangular
matrices in $GL_2(\cA)$ and $G_2^q \cong G_2^1 \cong (\cA,+)$ for all
$q \in \cA$.

For any $q \in \cA$, we get a bicrossed product l.c.\ quantum group
$(M,\de)_q$ whose right regular representation $W_q$ lives on the
constant Hilbert space $L^2(G_2^1 \times G_1^1)$. It is easy to check
that the application $q \mapsto W_q$ is strongly$^*$
continuous. So, we have a continuous family of l.c.\ quantum
groups. If $q \in \cAu$, it is clear that $(M,\de)_q \cong
(M,\de)_1$, because then, the first component of $\pi_q$ is an isomorphism. On the other hand, if $q = 0$, we observe that $G_2^0$ is
a normal subgroup of $G^0$. So, $G_1^0$ acts by automorphisms on
$(\cA,+)$ and $G^0 \cong G_1^0 \ltimes \cA$. So, $(M,\de)_0$ is
commutative and isomorphic to $(L^\infty(H),\de_H)$, where $H := G_1^0
\ltimes \hat{\cA}$ and $\hat{\cA}$ is the Pontryagin dual of
$(\cA,+)$.

In this precise sense, the l.c.\ quantum groups $(M,\de)_q$ are
continuous deformations of the l.c.\ group $H$.

If $\cA$ is commutative, we can remain closer to \cite{VV}, Section
5.4 and quotient out the center of $G^q$ consisting of the matrices
$\left(\begin{smallmatrix} a & 0 \\ 0 & a \end{smallmatrix}
\right)\subqq$, for $a \in \cAu$. 
\end{example}

\section{Concluding remarks} \label{sec.semireg}

\subsection{Pentagonal transformations} \label{subsec.pentagonal}

Following \cite{BS2}, we call $v : X \times X \recht X \times X$ a
pentagonal transformation, when $X$ is a (standard) measure
space and $v$ is a measure class isomorphism satisfying the pentagonal
relation $$v_{23} \na v_{13} \na v_{12} = v_{12} \na v_{23} \; .$$
Associated to any pentagonal transformation $v$, we have a
multiplicative unitary $V$ on the Hilbert space $L^2(X)$, defined by
$(V \xi)(x,y) = d(x,y)^{1/2} \xi(v(x,y))$, where $d$ is the right
Radon-Nikodym derivative.

Not all pentagonal transformations are nice: in
Remark~\ref{strange.mult.un}, we obtained a pentagonal transformation
such that $S_V$ and $\Sh_V$ are not even \cst-algebras.

In \cite{BS2}, there was given a natural sufficient condition for a
pentagonal transformation to be good. As explained in the
introduction, the result in \cite{BS2} is incorrect as stated, but can
be repaired as follows.

\begin{proposition} \label{prop.pentagonal}
Let $v$ be a pentagonal transformation on $X$.
Introduce binary composition laws on $X$ by writing $v(x,y) = (x
\bol y , x \dia y)$. Suppose that the transformations $\phi(x,y) = (x
\bol y, y)$ and $\eta(x,y) = (x,x \dia y)$ are measure class isomorphisms.

Then, there exists a matched pair $G_1,G_2$ of closed subgroups of a
l.c.\ group $G$ (in the sense of Definition~\ref{def.matched}),
commuting actions on $X$ of $G_2$ on the left and $G_1$ on the right,
and a $G_2 \times G_1$-equivariant measurable map $f : X \recht G$,
such that
$$v(x,y) = \bigl(x \cdot p_1\bigl(p_2(f(x))^{-1} f(y) \bigr),
p_2(f(x))^{-1} \cdot y\bigr) \qquad\text{for almost all}\quad (x,y) \in X
\times X \; .$$
\end{proposition}

We sketch the proof of this result, following \cite{BS2} (there only
appears an error in the proof of Proposition 3.4~a) of \cite{BS2}).  In order to find $G_1,G_2$ and
$G$, we write $v^{-1}(x,y) = (x \car y, x
\ster y)$ and further, $\psi'(x,y) = (y \ster x,y)$, $w(a,b,c,d) = (a
\bol (b \dia c) , d \ster (b \bol c) , c, d)$. One can check that
$\phi$ and $\psi'$ are pentagonal transformations on $X$ and $w$ is a
pentagonal transformation on $X \times X$. Because the pentagonal
transformations $\phi,\psi',w$ are all of the form $(r,s) \mapsto
(\ldots, s)$, it follows from Lemme 2.1 in \cite{BS2} that we can find
l.c.\ groups $G_1,G_2$ and $G$, right actions of $G_i$ on $X$ and of
$G$ on $X \times X$, and equivariant measurable maps $f_i : X \recht
G_i$, $F : X \times X \recht G$ such that
$$\phi(x,y) = (x \cdot f_1(y) , y) \; , \quad \psi'(x,y) = (x \cdot
f_2(y) , y) \quad\text{and}\quad w(a,b,c,d) = ((a,b) \cdot F(c,d),c,d)
\; .$$
One also verifies that $w = \psi'_{24} v_{21} \phi_{13} v_{21}^{-1}$.

Denote by $V$ the multiplicative unitary
associated to $v$ and by $V^1,V^2$ and $W$ the multiplicative
unitaries associated to $\phi,\psi'$ and $w$, resp.
So, we have $W=V^*_{21} V^1_{13} V_{21} V^2_{24}$. Hence, the group von Neumann
algebras of $G_1$ and $G_2$ are von Neumann subalgebras of the group
von Neumann algebra of $G$ and so, $G_1$ and $G_2$ can be considered
as closed subgroups of $G$. Under this identification, one has $F(c,d)
= f_1(c) f_2(d)$ almost everywhere. Because $F$ is surjective (from a
measure theoretic point of view), it follows that the
complement of $G_1 G_2$ has measure zero. From the formula for $W$, we
also get that $L^\infty(G) \subset L^\infty(G_1 \times G_2)$ and we
have to prove that we have an equality. Then, $G_1 \cap G_2 = \{e\}$
and we have a matched pair $G_1,G_2$ in $G$. Denote by $\weak$ the weak
closure. Then,
\begin{align*}
L^\infty(G) &= [(\om \ot \io \ot \io)(V^*_{21} V^1_{13} V_{21}
V^2_{24}) \mid \om \in \B(L^2(X\times X))_*]\weak \\ &= [(\om \ot \io \ot
\io)\bigl(V^1_{13} \bigl( (L^\infty(X) \ot 1)V(1 \ot L^\infty(X)) \bigr)_{21}
V^2_{24} \bigr) \mid \om \in \B(L^2(X\times X))_*]\weak \; .
\end{align*}
If $H,K \in L^\infty(X)$, we observe that $(H \ot 1) V (1 \ot K) V^*$
is the function $(x,y) \mapsto H(x)K(x \dia y)$. As the transformation
$(x,y) \mapsto (x,x \dia y)$ is a measure class isomorphism, the weak
closure of these functions is the whole of $L^\infty(X \times X)$. So,
\begin{align*}
L^\infty(G) &= [(\om \ot \io \ot \io)\bigl( V^1_{13} L^\infty(X)_1 V_{21} V^2_{24} \bigr) \mid \om \in \B(L^2(X\times
X))_*]\weak \\ &= [(\om \ot \io \ot \io)\bigl( V^1_{13} (\Sh_{V^1}'
L^\infty(X))_1 V_{21} V^2_{24} \bigr) \mid \om \in \B(L^2(X\times
X))_*]\weak \\ &= [ (L^\infty(G_1) \ot 1) (\om \ot \io \ot \io)(
V^2_{24} ) \mid \om \in \B(L^2(X\times X))_*]\weak = L^\infty(G_1
\times G_2) \; ,
\end{align*}
where we used that $\Sh_{V^1}' L^\infty(X)$ is weakly dense in $\B(L^2(X))$.
We conclude that we indeed get a matched pair $G_1,G_2$ in $G$.

Because $\eta'(x,y) = (y\dia x,y)$ is a measure class isomorphism, we
associate with it a unitary operator $B$ on $L^2(X \times X)$. One
checks that $\psi'_{23} \eta'_{12} \eta'_{13} = \eta'_{12} \psi'_{23}$
and hence, $B_{13} B_{12} V^2_{23} = V^2_{23} B_{12}$. It follows that
$B^*$ is a representation of $V^2$. Hence, $B \in \B(L^2(X)) \ot
L^\infty(G_2)$. It follows that there exists a left action $\actsl$ of
$G_2$ on $X$ such that $y \dia x = f_2(y) \actsl x$ almost everywhere.

Analogously, we can find a right action of $G_1$ on $X$ such that $x
\actsr f_1(y)^{-1} = x \car y$. The left action of $G_2$ and the right
action of $G_1$ commute and defining $f(x) = f_1(x) f_2(x)^{-1}$, we
have a $G_2 \times G_1$-equivariant measurable map $f : X \recht G$. One can show that
this map and the commuting actions of $G_1$ and $G_2$ satisfy the
conclusion of the proposition.

\begin{remark}
Combining Proposition~\ref{prop.pentagonal} and
Proposition~\ref{bicroscovariant}, we observe that the multiplicative
unitary $V$ associated to a \emph{good} pentagonal transformation $v$
is precisely a covariant image of the regular representation of the
associated bicrossed product.
In particular, $(S_V,\de_V)$ is always a bicrossed product l.c.\ quantum group.
\end{remark}

\begin{remark}
In the beginning of Section~\ref{sec.bicros}, we started with a
matched pair $G_1,G_2$ of closed subgroups of $G$ and constructed with
them a $^*$-automorphism $\tau : L^\infty(G_1 \times G_2) \recht
L^\infty(G_2 \times G_1)$ such that $\tau \si$ is a matching 
of $(L^\infty(G_2),\de_2)$ and 
$(L^\infty(G_1),\de_1^{\text{\rm op}})$ with trivial cocycles, in the
sense of \cite{VV}, Definition~2.1. Suppose now that, conversely,
$\tau$ is such that $\tau \si$ is a matching. Hence, $\tau$ is a
faithful $^*$-homomorphism satisfying
$$(\tau \ot \io)(\io \ot \tau)(\de_1 \ot \io) = (\io \ot \de_1) \tau
\quad\text{and}\quad (\io \ot \tau)(\tau \ot \io)(\io \ot \de_2) =
(\de_2 \ot \io)\tau \; .$$
Defining $\al$ and $\be$ to be the restrictions of $\tau$ to
$L^\infty(G_1)$ and $L^\infty(G_2)$ respectively, we get a left action
of $G_2$ on $L^\infty(G_1)$ and a right action of $G_1$ on
$L^\infty(G_2)$. Hence, we can write $\tau(F)(s,g) =
F(\al_s(g),\be_g(s))$ for almost all $g,s$ and where $(\al_s)$ and
$(\be_g)$ are actions of $G_2$, resp.\ $G_1$ by
measure class isomorphisms of $G_1$, resp. $G_2$. We can define
with the same formula as in Definition~3.3, the multiplicative unitary
$W$, which corresponds to the pentagonal transformation
$$v(s,g,t,h) = (s,g \; \al_{\be_g(s)^{-1} t}(h), \be_g(s)^{-1} t,h)$$
on $X:=G_2 \times G_1$. Corresponding to $v$, we have the mappings
$\phi$ and $\eta$ as above. It is clear that $\eta$ is a measure class
isomorphism. We prove that the same holds for $\phi$.
Because we have a pentagonal
transformation, we can define the (measure theoretically) surjective
mapping $P(x,y) = x \bol y$. Defining the measure class isomorphism
$u_\al(s,g) = (s,\al_s(g))$ and the surjective mapping $P_1(s,g,t,h) =
(s,gh)$, one checks that $u_\al P = P_1 (u_\al \times u_\al)$, where we
use the relation $\al_s(gh) = \al_s(g) \al_{\be_g(s)}(h)$ for almost
all $s,g,h$, which follows from $(\io \ot \de_1)\al = (\tau \ot
\io)(\io \ot \al)\de_1$. It follows that $P = u_\al^{-1} P_1 (u_\al
\times u_\al)$ and so $\phi = (u_\al^{-1} \times u_\al^{-1}) \phi^1_{24}
(u_\al \times u_\al)$, where $\phi^1(g,h) = (gh,h)$. Hence, also
$\phi$ is a measure class isomorphism and we can apply
Proposition~\ref{prop.pentagonal}. We find a l.c.\ group $G$ and
going through the construction, we
observe that we identify $G_1$ and $G_2$
with closed subgroups of $G$ such that $G_1,G_2$ is a matched pair and
$\al_s(g) = p_1(sg)$, $\be_g(s) = p_2(sg)$.
\end{remark}

\subsection{Regularity and semi-regularity}

We will give several characterizations of regularity and
semi-regularity of l.c.\ quantum groups. In fact, it will become clear
that it is very amazing that there indeed exist non-semi-regular l.c.\
quantum groups.

\begin{terminology}
A l.c.\ quantum group is said to be regular, resp.\ semi-regular, if
its right (or, equivalently, left) regular representation is a
regular, resp.\ semi-regular multiplicative unitary.
\end{terminology}

Let $V$ be the right regular representation of a fixed l.c.\ quantum
group $(M,\de)$, as in Section~\ref{sec.prelim}. Denote by $S$ its
underlying \cst-algebra as in Equation~\eqref{eq.SandSh}.

In the proof of
Proposition~\ref{char}, we introduced the modular conjugations $J$ and
$\Jh$, that we will use extensively. They are anti-unitary operators
on the Hilbert space $H$, which is the GNS-space of the invariant
weights. Observe that $J M J = M'$, $\Jh \Mh \Jh = {\Mh}'$, $\Jh S \Jh
= S$ and $J \Sh J = \Sh$.

When we write $\K$, we always mean $\K(H)$.

From Proposition~\ref{char}, we know that regularity is equivalent
with $\K = [S \; \Sh]$ and this is not always the case. However, if
we put slightly more, we do get $\K$.

\begin{lemma} \label{lemma.compact}
We have $[JSJ \; S \; \Sh]=\K$.
\end{lemma}
\begin{proof}
From Equation~\eqref{eq.irred} and using the notation $U=J \Jh$, we
know that, up to a scalar
$$(1 \ot U) \Si = (U \ot U) V^* (1 \ot U) \Si V^* \Si (U \ot 1) V^* \;
.$$
Applying $\om \ot \io$, it follows that
$$\K = [(\om \ot \io)\bigl( (U \ot U) V^* (1 \ot U) \Si V^* \Si (U
\ot 1) V^* \bigr) \mid \om \in \B(H)_* ] \; .$$
Because $\K = [J S J \; \K \; S]$ and because $V \in \M(\K \ot
S)$, it follows that
$$\K = [JSJ \; (\om \ot \io)(\Si V^* \Si) \; S \mid \om \in
\B(H)_*] = [JSJ \; \Sh \; S] \; .$$
Because $[\Sh \; S]$ is a \cst-algebra, we are done.
\end{proof}

Recall that $V \in \M(\K \ot S)$ and $S = [(\om \ot \io)(V) \mid \om
\in \B(H)_*]$. So, the following result is quite surprising.

\begin{proposition}
The l.c.\ quantum group $(M,\de)$ is semi-regular if and only if
$$[(\om \ot \io)(V^* (1 \ot S) V) \mid \om \in \B(H)_* ] \cap S
\neq \{0\} \; .$$
If this is the case, we have $S \subset [(\om \ot \io)(V^* (1 \ot S)
V) \mid \om \in \B(H)_* ]$. The quantum group is regular if and only
if the equality holds.
\end{proposition}
\begin{proof}
We first claim that in general
$$[\Sh \; JSJ \; (\om \ot \io)(V^* (1 \ot S) V) \mid \om \in \B(H)_* ]
= [\Sh \; JSJ] \; .$$
We make the following calculation, where we use the pentagonal
equation, Lemma~\ref{lemma.compact} and the formula $V^* = (J \ot
\Jh)V(J \ot \Jh)$.
\begin{align*}
[\Sh \; JSJ \; (\om \ot \io)(V^* (1 \ot S) V) \mid \om \in \B(H)_* ]
&= [(\om \ot \io \ot \mu)(V_{23} V_{12}^* (1 \ot JSJ \; S \ot 1)
V_{12} ) \mid \om,\mu \in \B(H)_* ] \\
&= [(\om \ot \io \ot \mu)(V_{13}^* V_{12}^* V_{23} (1 \ot JSJ \; S \ot 1)
V_{12} ) \mid \om,\mu \in \B(H)_* ] \\
&= [(\om \ot \io)(V^* (1 \ot \Sh \; JSJ \; S) V) \mid \om \in \B(H)_*
] \\ &= \Jh [(\io \ot \om)(\Si V (1 \ot \K) V^* \Si) \mid \om \in
\B(H)_*] \Jh \\ &= \Jh [\cC(V) \cC(V)^*] \Jh = [\Sh \; JSJ] \; ,
\end{align*}
where, in the end, we used the calculation of $\cC(V)$ from the proof
of Proposition~\ref{char}. This proves our claim.

Suppose that $a \neq 0$ and $a \in [(\om \ot \io)(V^* (1 \ot S) V)
\mid \om \in \B(H)_* ] \cap S$. Take $b \in S$ and $x \in \Sh$ such
that $x JbJ a \neq 0$. Combining our claim with
Lemma~\ref{lemma.compact}, we see that $x JbJ a \in [\Sh \; JSJ] \cap
\K$. Hence, $[\Sh \; JSJ] \cap
\K \neq \{0\}$. Because the \cst-algebra $[\Sh \; JSJ]$ acts
irreducibly on $H$, it follows that $\K \subset [\Sh \; JSJ]$ and
so, the l.c.\ quantum group is semi-regular.

Conversely, suppose that the l.c.\ quantum group is semi-regular.
Then,
\begin{align*}
[(\om \ot \io)(V^* (1 \ot S) V) \mid \om \in \B(H)_* ] &= [(\mu \ot
\om \ot \io)(V^*_{23} V^*_{13} V_{23}) \mid \om,\mu \in \B(H)_* ] \\ &= 
[(\mu \ot
\om \ot \io)(V^*_{23} V_{12} V_{23}) \mid \om,\mu \in \B(H)_* ] \\
&= [(\om \ot \io)(V^* (S \ot 1) V) \mid \om \in \B(H)_* ] \\ &= 
[(\om \ot \io)(V^* (\Jh \Sh \Jh \; S \ot 1) V) \mid \om \in \B(H)_* ]
\\ & \supset [(\om \ot \io)(V^* (\K \ot 1) V) \mid \om \in \B(H)_*
] = S \; .
\end{align*}
If the l.c.\ quantum group is regular, the $\supset$ above is an
equality. If, finally, $[(\om \ot \io)(V^* (1 \ot S) V) \mid \om \in
\B(H)_*] = S$, it follows from our claim above and Lemma~\ref{lemma.compact}
that
$$[\Sh \; JSJ] = [\Sh \; JSJ \; S] = \K \; .$$ This precisely gives
the regularity.
\end{proof}

Another characterization of semi-regularity is in terms of coactions
of l.c.\ quantum groups.

Let $B$ be a \cst-algebra and $\al : B \recht \M(S \ot B)$ a
non-degenerate $^*$-homomorphism satisfying $(\io \ot \al)\al = (\de
\ot \io)\al$. We say that $\al$ is a \emph{coaction} of $S$ on
$B$.

A seemingly natural definition for the \emph{continuous} elements
under the coaction is
$$T:=[(\om \ot \io)\al(B) \mid \om \in \B(H)_* ] \; .$$
\begin{proposition}
The l.c.\ quantum group $(M,\de)$ is semi-regular if and only if the
above closed linear space $T$ of continuous elements is a \cst-algebra
for all coactions $\al$ of $S$ on a \cst-algebra $B$.
\end{proposition}
\begin{proof}
Suppose first that $(M,\de)$ is semi-regular. Then,
\begin{align*}
T & \supset [(\om \ot \io)\al \bigl( y \; (\mu \ot \io)\al(x) \bigr) \mid
\om,\mu \in \B(H)_*, x,y \in B] \\ &= [(\mu \ot \om \ot
\io)(\al(y)_{23} \; V_{12} \al(x)_{13} V^*_{12}) \mid
\om,\mu \in \B(H)_*, x,y \in B] \\ &= [(\mu \ot \om \ot
\io)(\al(y)_{23} \bigl( (\K \ot 1) V (1 \ot \K) \bigr)_{12}
\al(x)_{13}) \mid \om,\mu \in \B(H)_*, x,y \in B] \\ & \supset 
[(\mu \ot \om \ot
\io)(\al(y)_{23} (\K \ot \K \ot 1)
\al(x)_{13} ) \mid \om,\mu \in \B(H)_*, x,y \in B] = [TT] \; .
\end{align*}
Hence, $T$ is an algebra. It is clear that $T^*=T$ and so, $T$ is a
\cst-algebra.

If conversely, the space of continuous elements is always a \cst-algebra,
we define $B:= \K(H \oplus \C)$. Consider the left regular
representation $W \in \M(S \ot \K)$ satisfying $(\de \ot \io)(W) =
W_{13} W_{23}$ and put $X:=W \oplus 1 \in \M(S \ot B)$. Define $\al(x)
= X^*(1 \ot x) X$ for all $x \in B$. Writing $\te_\xi$ for the obvious
operator in $\K(\C,H)$ whenever $\xi \in H$, we observe that
$$T = \left\{ \begin{pmatrix} [(\om \ot \io)(W^* (1 \ot \K) W) \mid \om \in
  \B(H)_*] & \te_\xi \\ \te_\eta^* & \C \end{pmatrix} \mid \xi,\eta
\in H \right\} \; .$$
So, because $T$ is a \cst-algebra, we get
$$\K \subset [(\om \ot \io)(W^* (1 \ot \K) W) \mid \om \in
  \B(H)_*] \; .$$
Using the fact that $W = \Si (U \ot 1) V (U^* \ot 1) \Si$, we conclude
that
$$\K \subset U [(\io \ot \om)(V^* \Si (1 \ot \K) \Si V) \mid \om
\in \B(H)_* ] U^* = U [\cC(V)^* \cC(V)] U^* = [S \; \Sh] \; .$$
So, $(M,\de)$ is semi-regular.
\end{proof}

We conclude that the notion of continuous elements is problematic for coactions of
non-semi-regular quantum groups. Hence, it is not surprising that is
equally problematic to define \emph{continuous} coactions.

There are at least two natural definitions. Let $\al : B \recht \M(S
\ot B)$ be a coaction.
\begin{itemize}
\item We call $\al$ continuous in the weak sense when $B=[(\om \ot \io)\al(B)
  \mid \om \in \B(H)_* ]$.
\item We call $\al$ continuous in the strong sense when $S \ot B =
  [\al(B)(S \ot 1)]$.
\end{itemize}

We have the following result.

\begin{proposition}
Consider a coaction $\al : B \recht \M(S \ot B)$ of $S$ on $B$.
If $(M,\de)$ is regular, continuity in the weak sense implies
continuity in the strong sense.

If $(M,\de)$ is semi-regular, but not regular, there exists a coaction
$\al$ that is continuous in the weak sense, but not in the strong
sense.
\end{proposition}
\begin{proof}
If $(M,\de)$ is regular, we know from \cite{BS} that $[(\K \ot 1)V(1
\ot S)] = \K \ot S$. So, if $\al$ is continuous in the weak sense, we
get
\begin{align*}
[\al(B) (S \ot 1)] &= [\al((\om \ot \io)\al(B)) (S \ot 1) \mid \om \in
\B(H)_*] \\ &= [(\om \ot \io \ot \io)(V_{12} \al(B)_{13} V^*_{12} (\K
\ot S \ot 1)) \mid \om \in \B(H)_*] \\ &= [(\om \ot \io \ot
\io)(V_{12} \al(B)_{13} (1 \ot S \ot 1)) \mid \om \in \B(H)_*] \\ &= [
(\om \ot \io \ot \io) \bigl( ((\K \ot 1)V(1 \ot S))_{12} \al(B)_{13}
\bigr) \mid \om \in \B(H)_*] \\ & = S \ot [(\om \ot \io)\al(B) \mid
\om \in \B(H)_*] = S \ot B \; .
\end{align*}
Hence, $\al$ is continuous in the strong sense.

Suppose now that $(M,\de)$ is semi-regular, but not regular. Define
$$B:= \left\{ \begin{pmatrix} [\Sh \; S] & \te_\xi \\ \te_\eta^* & \C
  \end{pmatrix} \mid \eta,\xi \in H \right\} \quad\text{with}\quad
\al \begin{pmatrix} x & \te_\xi \\ \te_\eta^* & \lambda \end{pmatrix}
= \begin{pmatrix} W^*(1 \ot x)W & W^*(1 \ot \te_\xi) \\
(1 \ot \te_\eta^*)W & 1 \ot \lambda \end{pmatrix} \; ,$$
where $W \in \M(S \ot \Jh \Sh \Jh)$ is again the left regular
representation. It is easy to verify that $\al$ is a coaction which is
continuous in the weak sense, but not in the strong sense.
\end{proof}

\end{document}